\documentclass{tac}

\usepackage{paralist}
\usepackage{mathtools}
\usepackage{amsmath}
\usepackage{amssymb}
\usepackage{pgf,tikz}
\usepackage{tikz-cd}
\usepackage{mathrsfs}
\usepackage{stmaryrd}
\usepackage{quiver}

\usepackage{xy}

\input diagxy

\usepackage[colorlinks=true]{hyperref}
\hypersetup{allcolors=[rgb]{0.1,0.1,0.4}}

\author{Jan Jurka}

\thanks{The author acknowledges the support of the Grant Agency of the Czech Republic under the grant 22-02964S and of Masaryk University under the grant MUNI/A/1457/2023.}

\address{Department of Mathematics and Statistics, Faculty of Science, Masaryk University}

\title[Enriched small object argument]{An enriched small object argument over a cofibrantly generated base}

\copyrightyear{2025}


\eaddress{jan.jurka@mail.muni.cz}

\newtheorem{thm}{Theorem}
 
\newtheorem{lem}{Lemma}

\newtheoremrm{rem}{Remark}

\let\pf\proof
\let\epf\endproof
\mathrmdef{Hom}
\mathbfdef{Set}

\begin{document}

\maketitle
\begin{abstract}
The small object argument is a method for transfinitely constructing weak factorization systems originally motivated by homotopy theory. We establish a variant of the small object argument that is enriched over a cofibrantly generated weak factorization system. This enriched variant of the small object argument subsumes the ordinary small object argument for categories and also certain variants of the small object argument for 2-categories, (2,1)-categories, dg-categories and simplicially enriched categories.
\end{abstract}


\section{Introduction}\label{sec-Introduction}

Quillen \cite{Qui67} introduced a way of transfinitely constructing weak factorization systems, dubbed the \emph{small object argument}. The original motivation for the argument comes from the theory of model categories, which is a categorical approach to homotopy theory. Moreover, later on various variants of the small object argument (\cite[1.37]{AR94}, \cite{Gar09}, \cite[13.2.1]{Rie14}) became an important tool in category theory itself and also in other areas of mathematics such as model theory due to the connection between the argument and ubiquitous notions of injectivity and orthogonality.\par
Enriched category theory is part of category theory that deals with ``categories'' in which hom-sets are not necessarily simply sets anymore, but instead they are objects in some monoidal category (a base of enrichment). It is the purpose of this paper to find a variant of the small object argument in the context of enriched category theory. The enrichment will be over a cofibrantly generated weak factorization system. As special cases we obtain the classical 1-categorical small object argument for weak factorization systems, the 1-categorical small object argument for orthogonal factorization systems, and certain variants of the small object argument for 2-categories, (2,1)-categories, dg-categories and simplicially enriched categories.\par
The basic idea of the construction is as follows. Recall that a morphism $f \colon A \to B$ is said to have the \emph{left lifting property} with respect to a morphism $k \colon C \to D$, which we denote $f \mathrel{\pitchfork} k$, if for each commutative square
\[\begin{tikzcd}
	A & C \\
	B & D
	\arrow["f"', from=1-1, to=2-1]
	\arrow["r", from=1-1, to=1-2]
	\arrow["s"', from=2-1, to=2-2]
	\arrow["k", from=1-2, to=2-2]
\end{tikzcd}\]
there exists a diagonal $d \colon B \to C$ making the two triangles below commute.
\[\begin{tikzcd}
	A & C \\
	B & D
	\arrow["f"', from=1-1, to=2-1]
	\arrow["r", from=1-1, to=1-2]
	\arrow["s"', from=2-1, to=2-2]
	\arrow["k", from=1-2, to=2-2]
	\arrow["d", from=2-1, to=1-2]
\end{tikzcd}\]
Given a $\mathcal{V}$-enriched category $\mathcal{K}$, there is an object $\mathrm{Sq}(f,k)$ in $\mathcal{V}$ (of ``squares connecting $f$ to $k$'') for each pair of morphisms $f \colon A \to B$, $k \colon C \to D$ from the underlying category $\mathcal{K}_0$, and this makes the category of morphisms in the underlying category of $\mathcal{K}$ a $\mathcal{V}$\nobreakdash-category with $\mathrm{Sq}(f,k)$ serving as a hom-object. The object $\mathrm{Sq}(f,k)$ can be viewed as a pullback of $\mathcal{K}(f,D)$, $\mathcal{K}(A,k)$, which leads to an induced map $\langle f,k \rangle \colon \mathcal{K}(B,C) \to \mathrm{Sq}(f,k)$. This map is an isomorphism if and only if $f$ has the left lifting property with respect to $k$ in the enriched sense (explicitly defined in \cite{LW14}, implicitly defined in \cite{Day74}). Generalizing this, given a class $\mathcal{J}$ of morphisms in $\mathcal{V}_0$, asking that each morphism from $\mathcal{J}$ has an ordinary left lifting property with respect to $\langle f,k \rangle$ encodes a form of an enriched left lifting property of $f$ with respect to $k$ that is relative to $\mathcal{J}$. Our variant of the small object argument then involves an enriched category $\mathcal{K}$, a base of enrichment $\mathcal{V}$, a class $\mathcal{J}$ of morphisms in $\mathcal{V}_0$, and a class $\mathcal{I}$ of morphisms in $\mathcal{K}_0$. In order to perform the small object argument, we find suitable conditions on $\mathcal{K}$, $\mathcal{V}$, $\mathcal{J}$, and $\mathcal{I}$, and show that the small object argument can be performed under these conditions.\par
In the ordinary 1-categorical small object argument, one constructs the first morphism in the factorization by taking a (nested) transfinite composite of pushouts of morphisms from $\mathcal{I}$, and the second morphism in the factorization by using the universal property of a transfinite composite. In our enriched variant of the small object argument we replace pushouts by ``copowered pushouts'', which are relative to each morphism in $\mathcal{J}$. As a consequence of that, in each step of the transfinite construction we use $|\mathcal{J}|$-many different kinds of ``pushouts'' and this is done by cycling through $\mathcal{J}$.\par
The paper is organized in the following way: In the second section we recall some preliminaries, in the third section we mention an associativity of certain copowers that we use in the small object argument, in the fourth section we apply the associativity from the previous section to the enriched category of arrows, in the fifth section we define the needed notions such as enriched liftings, in the sixth section we prove the required stability properties, in the seventh section we prove a factorization lemma in the category of morphisms, and in the eighth section we finally perform the small object argument.
\newline\newline
{\bf Acknowledgements.} I would like to thank my doctoral advisor John Bourke for suggesting to me an interesting topic to investigate, for many fruitful discussions, for carefully reading drafts of the paper, and for many helpful suggestions on how to improve presentation. I am also grateful to the anonymous referee for many suggestions which substantially improved the presentation. Furthermore, I would like to thank Simon Henry for a helpful discussion on the broad picture of enriched small object arguments, and Nathanael Arkor for telling me about locally graded categories.

\section{Preliminaries}\label{sec-preliminaries}

In this short section we recall some well-known preliminary notions.
\definition
For a class $\mathcal{J}$ of morphisms in a category, ${}^{\pitchfork}\mathcal{J}$ denotes the class of all morphisms that have the left lifting property with respect to all morphisms from $\mathcal{J}$ and $\mathcal{J}^\pitchfork$ denotes the class of all morphisms that have the right lifting property with respect to all morphisms from $\mathcal{J}$. For the definition of $\mathcal{J}^\pitchfork$ recall that a morphism $k$ is said to have the \emph{right lifting property} with respect to a morphism $f$ if $f \mathrel{\pitchfork} k$.
\enddefinition
\definition
A \emph{weak prefactorization system} $\mathscr{F} = (\mathscr{L},\mathscr{R})$ is a pair of classes of morphisms in a category such that $\mathscr{L}^\pitchfork = \mathscr{R}$ and $\mathscr{L} = {}^{\pitchfork}\mathscr{R}$.\par
Furthermore, a weak prefactorization system $\mathscr{F}$ is called \emph{a weak factorization system} if for each morphism $f$ there exists a pair of morphisms ${g \in \mathscr{L}}$, ${h \in \mathscr{R}}$ such that $f = h \cdot g$.\par
Moreover, a weak factorization system $\mathscr{F}$ is said to be \emph{cofibrantly generated} if there exists a set $\mathcal{J}$ of morphisms such that $\mathscr{R} = \mathcal{J}^\pitchfork$.
\enddefinition
\definition
Suppose that $\lambda$ is an ordinal and $\mathcal{C}$ is a category that has colimits that are indexed by ordinals $\gamma \leq \lambda$. Then a \emph{$\lambda$-sequence} $D \colon \lambda \to \mathcal{C}$ is a functor such that for all limit ordinals $\alpha < \lambda$ the induced morphism $\operatorname{colim}_{\beta < \alpha} D\beta \to D\alpha$ is an isomorphism. Furthermore, we call the colimit injection $D0 \to \operatorname{colim} D$ the \emph{transfinite composition of the $\lambda$-sequence $D$.}
\enddefinition
\rem
If $(\mathscr{L},\mathscr{R})$ is a weak prefactorization system, then it is well-known and easy to see that $\mathscr{R}$ is stable under pullbacks, transfinite cocomposites, and isomorphisms, and dually $\mathscr{L}$ is stable under pushouts, transfinite composites, and isomorphisms. In particular, both classes are stable under binary composites.
\endrem
\rem\label{rem:adj}
If a functor $F \colon \mathcal{C} \to \mathcal{D}$ between categories is left adjoint to a functor $G \colon \mathcal{D} \to \mathcal{C}$, then the following equivalence holds for all morphisms $f$ in $\mathcal{C}$, $g$ in $\mathcal{D}$: $$F(f) \mathrel{\pitchfork} g \text{ if and only if } f \mathrel{\pitchfork} G(g).$$
\endrem
\section{Copowers in categories of $\mathcal{V}$-functors}

Throughout the paper we will assume that $\mathcal{V} = (\mathcal{V}_0,\otimes,I)$ is a \emph{cosmos}, i.e.\ a bicomplete symmetric monoidal closed category, and that $\mathcal{K}$ is a $\mathcal{V}$-category. Furthermore, in this section we will assume that $\mathcal{K}$ is a copowered $\mathcal{V}$-category admitting coends of the form (\ref{def:action}). For objects $V$ in $\mathcal{V}$ and $K$ in $\mathcal{K}$ we will denote the copower by $V \odot K$. Recall that $\mathcal{V}$ is a copowered $\mathcal{V}$-category in which copowers $U \odot V$ are given by the monoidal product $U \otimes V$, and hence from now on we will use the copower notation instead of the monoidal notation.\par
Let $\mathcal{A}$ be a small monoidal $\mathcal{V}$-category. We are going to show that the underlying category $[\mathcal{A},\mathcal{K}]_0$ of the $\mathcal{V}$-category $[\mathcal{A},\mathcal{K}]$ of $\mathcal{V}$-functors $\mathcal{A} \to \mathcal{K}$ is a copowered $[\mathcal{A},\mathcal{V}]$\nobreakdash-category. This is a quite straightforward consequence of the theory of enriched actegories, and it is the purpose of this section to explain this. Given two $\mathcal{V}$-categories $\mathcal{L}$ and $\mathcal{L}'$ we will denote by $\mathcal{L} \otimes \mathcal{L}'$ their tensor product \cite[p.\ 12]{Kel05}. Given two $\mathcal{V}$-functors $F, G \colon \mathcal{A} \to \mathcal{V}$, denote by $F \odot G$ the composite 
\[\begin{tikzcd}
	{\mathcal{A} \otimes \mathcal{A}} & {\mathcal{V} \otimes \mathcal{V}} & {\mathcal{V}}
	\arrow["{F \otimes G}", from=1-1, to=1-2]
	\arrow["\odot", from=1-2, to=1-3]
\end{tikzcd}\]
and by $m \colon \mathcal{A} \otimes \mathcal{A} \to \mathcal{A}$ the monoidal product on $\mathcal{A}$. Recall that the Day convolution $F*G \colon \mathcal{A} \to \mathcal{V}$ is defined \cite[(3.1)]{Day70} by:
$$(F*G)(x) := \int^{a,b \in \mathcal{A}} \big(\mathcal{A}(m(a,b),x) \odot F(a)\big) \odot G(b),$$
and can be characterized \cite[Definition 21.4]{MMSS01} as the $\mathcal{V}$-functor part of the left Kan extension of $F \odot G$ along $m$. Furthermore, recall that the Day convolution is a monoidal product on the $\mathcal{V}$-category $[\mathcal{A},\mathcal{V}]$, which follows from \cite[Theorem 3.3]{Day70}. Given a $\mathcal{V}$-functor $X \colon \mathcal{A} \to \mathcal{K}$, we define the $\mathcal{V}$-functor $F*X \colon \mathcal{A} \to \mathcal{K}$ by:
\begin{equation}\label{def:action}
(F * X)(x) := \int^{a,b \in \mathcal{A}} \big(\mathcal{A}(m(a,b),x) \odot F(a)\big) \odot X(b),
\end{equation}
and again $F * X$ can be equivalently characterized as the $\mathcal{V}$-functor part of the left Kan extension of
\[\begin{tikzcd}
	{F \odot X \colon\mathcal{A} \otimes \mathcal{A}} & {\mathcal{V} \otimes \mathcal{K}} & {\mathcal{K}}
	\arrow["{F \otimes X}", from=1-1, to=1-2]
	\arrow["\odot", from=1-2, to=1-3]
\end{tikzcd}\]
along $m$. Note that this is a notion analogous to the Day convolution where we now use copowers in $\mathcal{K}$ instead of the monoidal product in $\mathcal{V}$.
\definition
Suppose that $(\mathcal{M},*_\mathcal{M},I)$ is a monoidal $\mathcal{V}$-category. A \emph{left $\mathcal{M}$-actegory} is a $\mathcal{V}$-category $\mathcal{L}$ equipped with a $\mathcal{V}$-functor $* \colon \mathcal{M} \otimes \mathcal{L} \to \mathcal{L}$ together with $\mathcal{V}$-natural isomorphisms 
\begin{align*}
\alpha_{M,N,L} \colon M*(N*L) &\xrightarrow{\cong} (M*_{\mathcal{M}}N)*L, \\ 
\lambda_{L} \colon I * L &\xrightarrow{\cong} L
\end{align*}
satisfying coherence conditions \cite[(1.1), (1.2), (1.3)]{JK01}.
\enddefinition
\proposition\label{prop:act}
Suppose that $\mathcal{V}$ is a cosmos and $\mathcal{K}$ is a copowered $\mathcal{V}$-category admitting coends of the form (\ref{def:action}). Then 
\begin{enumerate}[(i)]
\item the $\mathcal{V}$-category $[\mathcal{A},\mathcal{K}]$ is a left $[\mathcal{A},\mathcal{V}]$-actegory such that the action on a fixed object of $[\mathcal{A},\mathcal{K}]$ always has a right adjoint, and
\item the category $[\mathcal{A},\mathcal{K}]_0$ is a copowered $[\mathcal{A},\mathcal{V}]$-category.
\end{enumerate}
\endproposition
\pf
In the proof we omit the verification of coherence conditions, since it is not needed for the purposes of our small object argument. We begin by proving the first part of the theorem, the second part will then be a corollary of the first part.\par
Let $F, G \colon \mathcal{A} \to \mathcal{V}$, $X \colon \mathcal{A} \to \mathcal{K}$ be $\mathcal{V}$-functors. The action on $[\mathcal{A},\mathcal{K}]$ that we are looking for is $*$ from Definition (\ref{def:action}). First we will show associativity of the action, i.e.\ that
\begin{equation}\label{eqn:astassoc}(F*G)*X \cong F*(G*X).\end{equation}
We will decorate isomorphisms that follow from the Yoneda Lemma by $\mathfrak{Y}$. We have the following chain of isomorphisms:
\begin{small}
\begin{align*}
\big((F*G)*X\big)(x) &= \int^{a,b \in \mathcal{A}} \big(\mathcal{A}(m(a,b),x) \odot (F*G)(a)\big) \odot X(b)\\
&= \int^{a,b \in \mathcal{A}} \Big(\mathcal{A}(m(a,b),x) \odot \Big(\int^{c,d \in \mathcal{A}} \big(\mathcal{A}(m(c,d),a) \odot F(c)\big) \odot G(d)\Big)\Big) \odot X(b) \\
&\cong \int^{a, b, c, d \in \mathcal{A}} \Big(\mathcal{A}(m(a,b),x) \odot \big((\mathcal{A}(m(c,d),a) \odot F(c)) \odot G(d)\big)\Big) \odot X(b) \\
&\cong \int^{a, b, c, d \in \mathcal{A}} \Big(\big((\mathcal{A}(m(a,b),x) \odot \mathcal{A}(m(c,d),a)) \odot F(c)\big) \odot G(d)\Big) \odot X(b) \\
&\stackrel{\mathfrak{Y}}{\cong} \int^{b, c, d \in \mathcal{A}} \Big(\big(\mathcal{A}(m(m(c,d),b),x) \odot F(c)\big) \odot G(d)\Big) \odot X(b) \\
&\cong \int^{b, c, d \in \mathcal{A}} \Big(\big(\mathcal{A}(m(c,m(d,b)),x) \odot F(c)\big) \odot G(d)\Big) \odot X(b) \\
&\stackrel{\mathfrak{Y}}{\cong} \int^{a, b, c, d \in \mathcal{A}} \Big(\big((\mathcal{A}(m(c,a),x) \odot \mathcal{A}(m(d,b),a)) \odot F(c)\big) \odot G(d)\Big) \odot X(b) \\
&\cong \int^{a, b, c, d \in \mathcal{A}} \big(\mathcal{A}(m(c,a),x) \odot F(c)\big) \odot \Big(\big(\mathcal{A}(m(d,b),a) \odot G(d)\big) \odot X(b)\Big) \\
&\cong \int^{c,a \in \mathcal{A}} \big(\mathcal{A}(m(c,a),x) \odot F(c)\big) \odot \Big(\int^{d,b \in \mathcal{A}} \big(\mathcal{A}(m(d,b),a) \odot G(d)\big) \odot X(b)\Big)\\
&= \int^{c,a \in \mathcal{A}} \big(\mathcal{A}(m(c,a),x) \odot F(c)\big) \odot (G*X)(a)\\
&= \big(F*(G*X)\big)(x).
\end{align*}
\end{small}\par
Note that the eighth line in the chain above uses the symmetry of the monoidal product in $\mathcal{V}$.\par
Denote by $i$ the unit object of $\mathcal{A}$. We will show that $\mathcal{A}(i,-) \colon \mathcal{A} \to \mathcal{V}$ is the unit of the action $*$. This can be done in an analogous way as showing associativity by using the coend definition of the action $*$, however we will show it by using the Kan extension characterization of the action $*$:
\begin{align*}
[\mathcal{A}, \mathcal{K}]\big(\mathcal{A}(i,-) * X, Y\big) &\cong [\mathcal{A} \otimes \mathcal{A}, \mathcal{K}]\big(\mathcal{A}(i,-) \odot X, Y \cdot m\big) \\
&\cong \int_{a,b \in \mathcal{A}} \mathcal{K}\big(\mathcal{A}(i,a) \odot X(b), Y(m(a,b))\big) \\
&\cong \int_{a,b \in \mathcal{A}} \mathcal{V}\Big(\mathcal{A}(i,a), \mathcal{K}\big(X(b), Y(m(a,b))\big)\Big) \\
&\stackrel{\mathfrak{Y}}{\cong} \int_{b \in \mathcal{A}} \mathcal{K}\big(X(b), Y(m(i,b))\big) \\
&\cong \int_{b \in \mathcal{A}} \mathcal{K}\big(X(b), Y(b)\big) \\
&\cong [\mathcal{A}, \mathcal{K}](X, Y),
\end{align*}
and thus $\mathcal{A}(i,-) * X \cong X$ by the Yoneda Lemma. If $X, Y \colon \mathcal{A} \to \mathcal{K}$ are $\mathcal{V}$-functors, define
\begin{equation}\label{def:anglehom}\langle X,Y \rangle := \int_{a \in \mathcal{A}} \mathcal{K}\big(X(a), Y(m(-,a))\big).\end{equation}\par
In order to finish the proof, we will show that $\langle X, - \rangle \colon [\mathcal{A},\mathcal{K}] \to [\mathcal{A}, \mathcal{V}]$ is a right adjoint to $- * X \colon [\mathcal{A},\mathcal{V}] \to [\mathcal{A}, \mathcal{K}]$. Indeed,
\begin{align*}
[\mathcal{A}, \mathcal{K}](F * X, Y) &= \int_{c \in \mathcal{A}} \mathcal{K}\big((F*X)(c), Y(c)\big)\\
&\cong \int_{c \in \mathcal{A}} \mathcal{K}\Big(\int^{a,b \in \mathcal{A}} \big(\mathcal{A}(m(a,b),c) \odot F(a)\big) \odot X(b), Y(c)\Big)\\
&\cong \int_{a, b, c \in \mathcal{A}} \mathcal{V}\Big(\mathcal{A}(m(a,b),c) \odot F(a), \mathcal{K}(X(b),Y(c))\Big) \\
&\cong \int_{a, b, c \in \mathcal{A}} \mathcal{V}\Big(F(a), \mathcal{V}\big(\mathcal{A}(m(a,b),c),\mathcal{K}(X(b),Y(c))\big)\Big) \\
&\cong \int_{a \in \mathcal{A}} \mathcal{V}\Big(F(a), \int_{b,c \in \mathcal{A}} \mathcal{V}\big(\mathcal{A}(m(a,b),c),\mathcal{K}(X(b),Y(c))\big)\Big) \\
&\stackrel{\mathfrak{Y}}{\cong} \int_{a \in \mathcal{A}} \mathcal{V}\Big(F(a), \int_{b \in \mathcal{A}} \mathcal{K}\big(X(b), Y(m(a,b))\big)\Big) \\
&= \int_{a \in \mathcal{A}} \mathcal{V}\big(F(a), \langle X,Y \rangle (a)\big) \\
&\cong [\mathcal{A},\mathcal{V}]\big(F, \langle X,Y \rangle\big).
\end{align*}
By \cite[6.\ Appendix on tensored $\mathcal{V}$-categories]{JK01}, the second part of the theorem follows from the first part. We remark that copowers are given by the action $*$ from Definition (\ref{def:action}) and the hom-object for $X, Y \colon \mathcal{A} \to \mathcal{K}$ is given by $\langle X, Y \rangle$ from Definition (\ref{def:anglehom}).
\epf
\rem
When $\mathcal{V} = \mathbf{Set}$ with the cartesian monoidal structure, the enrichment of $[\mathcal{A},\mathcal{K}]_0$ over $[\mathcal{A},\mathbf{Set}]$ has been described by McDermott and Uustalu \cite{MU22}. They describe the enrichment directly in \cite[Definition 10]{MU22} using the language of locally $\mathcal{A}$-graded categories, which are an elementary formulation of $[\mathcal{A},\mathbf{Set}]$-enriched categories due to Wood \cite[Theorem 1.6]{Woo76}. For our purposes, the case of general $\mathcal{V}$ and the copowers are essential.
\endrem

\section{The enriched category of arrows}

We now specialise the results from the previous section to $\mathcal{A} = \mathbf{2}$, where $\mathbf{2}$ is the free $\mathcal{V}$-category on the category with two objects $0, 1$ and a single non-identity morphism $0 \to 1$. Recall that this free $\mathcal{V}$-category has as hom-objects the initial object and the unit object if the corresponding hom-sets are the empty set and the singleton, respectively. Furthermore, $\mathbf{2}$ can be equipped with the cartesian monoidal product $m$, which is given on objects $x, y \in \mathbf{2}$ by the formula $m(x,y) := \mathrm{min}(x,y)$.\par
As usual, we identify $\mathcal{V}_0$ with the underlying category of the $\mathcal{V}$-category $\mathcal{V}$. Then $[\mathbf{2},\mathcal{V}]$ and $[\mathbf{2},\mathcal{K}]$ are the $\mathcal{V}$-categories of morphisms in $\mathcal{V}_0$ and $\mathcal{K}_0$, respectively.\par
The fact that the category $[\mathbf{2},\mathcal{K}]_0$ is a copowered $[\mathbf{2},\mathcal{V}]$-category (i.e.\ part (ii) of Proposition \ref{prop:act} for $\mathcal{A} = \mathbf{2}$) and the explicit description of copowers and hom-objects appears in\ \cite[4.1.5]{Sub21} (note that this reference includes a local presentability assumption, which we do not assume).\par
Hom-objects in $[\mathbf{2},\mathcal{K}]$ as a $\mathcal{V}$-category are given by $\mathrm{Sq}(f,k)$ from the following definition.
\definition\label{def:sq}
For each pair of morphisms ${f \colon A \to B}$, ${k \colon C \to D}$ in $\mathcal{K}_0$, define $\mathrm{Sq}(f,k)$ to be the \emph{object of squares connecting $f$ to $k$}, i.e.\ the pullback-object in the following pullback square.
\[\begin{tikzcd}
	{\mathrm{Sq}(f,k)} & {\mathcal{K}(A,C)} \\
	{\mathcal{K}(B,D)} & {\mathcal{K}(A,D)}
	\arrow["{p_2}"', from=1-1, to=2-1]
	\arrow["{\mathcal{K}(f,D)}"', from=2-1, to=2-2]
	\arrow["{p_1}", from=1-1, to=1-2]
	\arrow["{\mathcal{K}(A,k)}", from=1-2, to=2-2]
\end{tikzcd}\]
Note that since $[\mathbf{2},\mathcal{K}]$ is a $\mathcal{V}$-category, we have the associated $\mathcal{V}_0$-valued hom-functor ${\operatorname{Sq}(-,-) \colon [\mathbf{2},\mathcal{K}]_0^{\mathrm{op}} \times [\mathbf{2},\mathcal{K}]_0 \to \mathcal{V}_0}$, and so in particular $\operatorname{Sq}(-,k) \colon [\mathbf{2},\mathcal{K}]_0^{\mathrm{op}} \to \mathcal{V}_0$ for a morphism $k$ in $\mathcal{K}_0$. We will make heavy use of this in what follows, and so now record its explicit description. If
\[\begin{tikzcd}
	A & K \\
	B & L
	\arrow["g", from=1-1, to=1-2]
	\arrow["{g'}"', from=2-1, to=2-2]
	\arrow["f"', from=1-1, to=2-1]
	\arrow["{f'}", from=1-2, to=2-2]
\end{tikzcd}\]
is a commutative square, i.e.\ a morphism $(g,g') \colon f \to f'$ in $[\mathbf{2},\mathcal{K}]_0$, then there is a unique morphism $\mathrm{Sq}((g,g'),k) \colon \mathrm{Sq}(f',k) \to \mathrm{Sq}(f,k)$ in $\mathcal{V}_0$ that makes the two top squares in the following diagram commute, since $(p_1,p_2)$ is a pullback.
\[\begin{tikzcd}
	{\mathcal{K}(K,C)} & {\mathrm{Sq}(f',k)} & {\mathcal{K}(L,D)} \\
	{\mathcal{K}(A,C)} & {\mathrm{Sq}(f,k)} & {\mathcal{K}(B,D)} \\
	& {\mathcal{K}(A,D)}
	\arrow["{\mathrm{Sq}((g,g'),k)}", dashed, from=1-2, to=2-2]
	\arrow["{p_2}"', from=2-2, to=2-3]
	\arrow["{p_1}", from=2-2, to=2-1]
	\arrow["{\mathcal{K}(g,C)}"', from=1-1, to=2-1]
	\arrow["{p'_1}"', from=1-2, to=1-1]
	\arrow["{p'_2}", from=1-2, to=1-3]
	\arrow["{\mathcal{K}(g',D)}", from=1-3, to=2-3]
	\arrow["{\mathcal{K}(A,k)}"', from=2-1, to=3-2]
	\arrow["{\mathcal{K}(f,D)}", from=2-3, to=3-2]
\end{tikzcd}\]
\enddefinition
\examples
\begin{enumerate}[(1)]
\item In the case $\mathcal{V} = \mathbf{Set}$, where $\mathbf{Set}$ is the monoidal category of sets in which the monoidal structure is given by the cartesian product, we get that $\mathcal{K}$ is a category. The elements of $\mathrm{Sq}(f,k)$ are commutative squares in $\mathcal{K}$ of the form:
\begin{equation}\label{cmd:fksquare}
\begin{tikzcd}
	A & C \\
	B & D
	\arrow["f"', from=1-1, to=2-1]
	\arrow["k", from=1-2, to=2-2]
	\arrow["r", from=1-1, to=1-2]
	\arrow["s"', from=2-1, to=2-2]
\end{tikzcd}
\end{equation}
We will call each such square a \emph{commutative $(f,k)$-square}.
\item Suppose that $\mathcal{V} = \mathbf{Cat}$, where $\mathbf{Cat}$ is the monoidal category of categories in which the monoidal structure is given by the cartesian product. This means that $\mathcal{K}$ is a 2-category. Then the objects of $\mathrm{Sq}(f,k)$ are commutative $(f,k)$-squares, and the morphisms in $\mathrm{Sq}(f,k)$ are pairs $\theta \colon r \Rightarrow r'$, $\theta' \colon s \Rightarrow s'$ of 2-cells such that $k \ast \theta = \theta' \ast f$.
\item If $\mathcal{V} = \mathbf{Grpd}$, where $\mathbf{Grpd}$ is the monoidal category of groupoids in which the monoidal structure is given by the cartesian product, then that means that $\mathcal{K}$ is a (2,1)-category. The objects of $\mathrm{Sq}(f,k)$ are commutative $(f,k)$-squares, and the morphisms in $\mathrm{Sq}(f,k)$ are pairs $\theta \colon r \Rightarrow r'$, $\theta' \colon s \Rightarrow s'$ of invertible 2-cells such that $k \ast \theta = \theta' \ast f$.
\item Suppose that $\mathcal{V} = \mathbf{Ch}$, where $\mathbf{Ch}$ is the monoidal category of chain complexes of left $R$-modules over a ring $R$ in which the monoidal structure is given by the tensor product of chain complexes. This means that $\mathcal{K}$ is a dg-category and $\mathrm{Sq}(f,k)$ is a chain complex
\[\begin{tikzcd}
	\cdots & {\mathrm{Sq}(f,k)_{n + 1}} & {\mathrm{Sq}(f,k)_n} & {\mathrm{Sq}(f,k)_{n -1}} & \cdots
	\arrow["{\partial_{n + 2}}", from=1-1, to=1-2]
	\arrow["{\partial_{n + 1}}", from=1-2, to=1-3]
	\arrow["{\partial_n}", from=1-3, to=1-4]
	\arrow["{\partial_{n - 1}}", from=1-4, to=1-5]
\end{tikzcd}\]
whose $n$-th degree elements are pairs $(r, s) \in \mathcal{K}(A,C)_n \times \mathcal{K}(B,D)_n$ that satisfy ${k \cdot r = s \cdot f \in \mathcal{K}(A,D)_n}$. The equality ${\partial_n(r,s) = (\partial_n(r), \partial_n(s))}$ defines differentials on $\mathrm{Sq}(f,k)$.
\item In the case $\mathcal{V} = \mathbf{SSet}$, where $\mathbf{SSet}$ is the monoidal category of simplicial sets with monoidal structure given by the pointwise cartesian product, we get that $\mathcal{K}$ is a simplicially enriched category and $\mathrm{Sq}(f,k)$ is a simplicial set whose $n$-simplices are pairs $(r, t) \in \mathcal{K}(A,C)_n \times \mathcal{K}(B,D)_n$ such that ${k \cdot r = t \cdot f \in \mathcal{K}(A,D)_n}$. Face maps and degeneracy maps on $\mathrm{Sq}(f,k)$ are defined by $d_{i}(r,t) = (d_{i}(r), d_{i}(t))$ and ${s_{i}(r,t) = (s_{i}(r), s_{i}(t))}$, respectively.
\end{enumerate}
\endexamples
In the remainder of this section we will assume that $\mathcal{K}$ is a copowered $\mathcal{V}$-category admitting pushouts of the form $\mathrm{dom}(u \mathrel{\Box} f)$ from the following definition.
\definition\label{def:uodotf}
Suppose that $u \colon U \to V$ is a morphism in $\mathcal{V}_0$ and $f \colon A \to B$ is a morphism in $\mathcal{K}_0$. Then $u \mathrel{\Box} f \colon \mathrm{dom}(u \mathrel{\Box} f) \to V \odot B$ is the induced morphism depicted in the following diagram, where $(i_1,i_2)$ is a pushout of $U \odot f$ and $u \odot A$.
\[\begin{tikzcd}
	{U \odot A} && {U \odot B} \\
	& {\mathrm{dom}(u \mathrel{\Box} f)} \\
	{V \odot A} && {V \odot B}
	\arrow["{u \mathrel{\Box} f}"', dashed, from=2-2, to=3-3]
	\arrow["{U \odot f}", from=1-1, to=1-3]
	\arrow["{u \odot B}", from=1-3, to=3-3]
	\arrow["{V \odot f}"', from=3-1, to=3-3]
	\arrow["{u \odot A}"', from=1-1, to=3-1]
	\arrow["{i_2}", from=3-1, to=2-2]
	\arrow["{i_1}"', from=1-3, to=2-2]
\end{tikzcd}\]
\enddefinition
\rem
Proposition \ref{prop:act} gives us that the category $[\mathbf{2},\mathcal{K}]_0$ of morphisms in $\mathcal{K}_0$ is a copowered $[\mathbf{2},\mathcal{V}]$-category whose copower action $*$ is given by $\Box$ from Definition \ref{def:uodotf}.
\endrem
\rem
The existence of coends of the form (\ref{def:action}) is equivalent to the existence of pushouts of the form $\mathrm{dom}(u \mathrel{\Box} f)$ from Definition \ref{def:uodotf}, since for $\mathcal{A} = \mathbf{2}$ these coends are of the form (\ref{def:action2}). The fact that $\Box$ is the copower action $*$ follows from the formula \begin{equation}\label{def:action2}(u * f)(x) = \int^{a,b \in \mathbf{2}} \big(\mathbf{2}(m(a,b),x) \odot u(a)\big) \odot f(b).\end{equation}
\endrem
\rem\label{rem:Vnab}
Since $\mathcal{V}$ is a copowered $\mathcal{V}$-category, it makes sense to also consider $\mathrm{Sq}(u,v)$ and $u \mathrel{\Box} v$ for morphisms $u, v$ in $\mathcal{V}_0$.
\endrem
Now from Equation (\ref{eqn:astassoc}) we obtain \begin{equation}\label{eqn:nabassoc}v \mathrel{\Box} (u \mathrel{\Box} f) \cong (v \mathrel{\Box} u) \mathrel{\Box} f,\end{equation} which will be useful for the purposes of our small object argument. Moreover, the hom\nobreakdash-object $\langle f, k \rangle$ for $[\mathbf{2},\mathcal{K}]$ as a $[\mathbf{2},\mathcal{V}]$-enriched category becomes $\langle f,k \rangle$ from the following definition, which will later be used to express enriched lifting properties.
\definition\label{def:inducede}
Suppose that $f \colon A \to B$, $k \colon C \to D$ are morphisms in $\mathcal{K}_0$. Then define $\langle f,k \rangle \colon \mathcal{K}(B,C) \to \mathrm{Sq}(f,k)$ to be the induced morphism depicted below.
\begin{equation}\label{cmd:inducede}
\begin{tikzcd}
	{\mathcal{K}(B,C)} \\
	& {\mathrm{Sq}(f,k)} & {\mathcal{K}(A,C)} \\
	& {\mathcal{K}(B,D)} & {\mathcal{K}(A,D)}
	\arrow["{p_1}", from=2-2, to=2-3]
	\arrow["{p_2}"', from=2-2, to=3-2]
	\arrow["{\mathcal{K}(A,k)}", from=2-3, to=3-3]
	\arrow["{\mathcal{K}(f,D)}"', from=3-2, to=3-3]
	\arrow["{\mathcal{K}(B,k)}"', from=1-1, to=3-2]
	\arrow["{\mathcal{K}(f,C)}", from=1-1, to=2-3]
	\arrow["{\langle f,k \rangle}"{description}, dashed, from=1-1, to=2-2]
\end{tikzcd}
\end{equation}
\enddefinition
\rem\label{rem:copbij}
We will denote by $-^*$ the bijective correspondence assigning to each morphism $U \to \mathcal{K}(A,B)$ in $\mathcal{V}_0$ a morphism $U \odot A \to B$ in $\mathcal{K}_0$, and by $-_*$ its inverse.
\endrem
\rem\label{rem:enabadj}
The fact that $\langle f,k \rangle$ is the hom-object for $[\mathbf{2},\mathcal{K}]$ as a $[\mathbf{2},\mathcal{V}]$-enriched category follows from the functor ${- \mathrel{\Box} f \colon [\mathbf{2},\mathcal{V}]_0 \to [\mathbf{2},\mathcal{K}]_0}$ being left adjoint to the functor $\langle f,- \rangle \colon [\mathbf{2},\mathcal{K}]_0 \to [\mathbf{2},\mathcal{V}]_0$ and the uniqueness of right adjoints. The adjunction is given as follows, in both cases given a commutative square on the left we obtain the commutative square on the right:
\[\begin{tikzcd}[column sep=3.2em]
	U & {\mathcal{K}(B,C)} && {\mathrm{dom}(u \mathrel{\Box} f)} & C \\
	V & {\mathrm{Sq}(f,k)} && {V \odot B} & D
	\arrow["{\langle f,k \rangle}", from=1-2, to=2-2]
	\arrow["v"', from=2-1, to=2-2]
	\arrow["w", from=1-1, to=1-2]
	\arrow["u"', from=1-1, to=2-1]
	\arrow["{u \mathrel{\Box} f}"', from=1-4, to=2-4]
	\arrow["k", from=1-5, to=2-5]
	\arrow["{(p_2 \cdot v)^*}"', from=2-4, to=2-5]
	\arrow["{(w^*, (p_1 \cdot v)^*)}", from=1-4, to=1-5]
\end{tikzcd}\]
and 
\[\begin{tikzcd}[column sep=3.2em]
	{\mathrm{dom}(u \mathrel{\Box} f)} & C && U & {\mathcal{K}(B,C)} \\
	{V \odot B} & D && V & {\mathrm{Sq}(f,k)}
	\arrow["{\langle f,k \rangle}", from=1-5, to=2-5]
	\arrow["{((g \cdot i_2)_*, h_*)}"', from=2-4, to=2-5]
	\arrow["{(g \cdot i_1)_*}", from=1-4, to=1-5]
	\arrow["u"', from=1-4, to=2-4]
	\arrow["{u \mathrel{\Box} f}"', from=1-1, to=2-1]
	\arrow["k", from=1-2, to=2-2]
	\arrow["h"', from=2-1, to=2-2]
	\arrow["g", from=1-1, to=1-2]
\end{tikzcd}\]
\endrem

\section{Enriched lifting properties}

From now on in the paper we will assume that the cosmos $\mathcal{V}$ is equipped with a weak prefactorization system $\mathscr{F} = (\mathscr{L}, \mathscr{R})$ on $\mathcal{V}_0$. The purpose of this section is to define basic notions needed for our small object argument: Enriched lifting properties, enriched weak factorization systems, and the notion of a class of morphisms in $\mathcal{V}_0$ being stable under corners.
\definition
Suppose that $f \colon A \to B$, $k \colon C \to D$ are morphisms in $\mathcal{K}_0$. Then we write $f \overset{\mathscr{F}}{\mathrel{\pitchfork}} k$ if $\langle f,k \rangle \in \mathscr{R}$, where $\langle f,k \rangle \colon \mathcal{K}(B,C) \to \mathrm{Sq}(f,k)$ is the induced morphism from Definition \ref{def:inducede}. \par Furthermore, if $f \overset{\mathscr{F}}{\mathrel{\pitchfork}} k$, then we say that $f$ has the \emph{left $\mathscr{F}$-lifting property} with respect to $k$, or equivalently that $k$ has the \emph{right $\mathscr{F}$-lifting property} with respect to $f$. Moreover, if $\mathcal{I}$ is a class of morphisms in $\mathcal{K}_0$, then we define
\begin{align*}
 \mathcal{I}^{\overset{\mathscr{F}}{\pitchfork}} := \{k \in [\mathbf{2},\mathcal{K}] \mid \forall f \in \mathcal{I} \colon f \overset{\mathscr{F}}{\mathrel{\pitchfork}} k\},\\
{}^{\overset{\mathscr{F}}{\pitchfork}}\mathcal{I} := \{f \in [\mathbf{2},\mathcal{K}] \mid \forall k \in \mathcal{I} \colon f \overset{\mathscr{F}}{\mathrel{\pitchfork}} k\}.
\end{align*}
\enddefinition
\examples\label{ex:liftings}
\begin{enumerate}[(1)]
\item\label{ex:liftingsSet} In the case $\mathcal{V} = \mathbf{Set}$, $\mathscr{F} = (\text{injective}, \text{surjective})$, we capture the ordinary (weak) lifting property, i.e.\ $f \overset{\mathscr{F}}{\mathrel{\pitchfork}} k$ iff $f \mathrel{\pitchfork} k$. Note that $\mathscr{F}$ is cofibrantly generated by $\mathcal{J} = \{u \colon \emptyset \to 1\}$. \par
If we instead choose the following ${\mathcal{J} = \{u \colon \emptyset \to 1, v \colon 2 \to 1\}}$ we capture the strong lifting property in which the diagonal is required to be unique because $\mathscr{F}$ becomes $(\text{all functions}, \text{bijections})$.
\item\label{ex:liftingsCat} Suppose that $\mathcal{V} = \mathbf{Cat}$, $\mathscr{F} = (\text{injective on objects}, \text{surjective equivalences})$. Then $f \overset{\mathscr{F}}{\mathrel{\pitchfork}} k$ iff $\langle f,k \rangle$ is a surjective equivalence, which happens iff for each pair of 1-cells $r \colon A \to C$, $s \colon B \to D$ satisfying $k \cdot r = s \cdot f$ there exists a diagonal $d \colon B \to C$ such that $d \cdot f = r$, $k \cdot d = s$, and furthermore if $d, d' \colon B \to C$ are 1-cells, and $\theta \colon d \cdot f \Rightarrow d' \cdot f$, $\theta' \colon k \cdot d \Rightarrow k \cdot d'$ are 2-cells such that $k \ast \theta = \theta' \ast f$, then there exists a unique 2-cell $\varphi \colon d \Rightarrow d'$ such that $\varphi \ast f = \theta$ and $k \ast \varphi = \theta'$. Note that $\mathscr{F}$ is cofibrantly generated by $\mathcal{J} = \{u \colon \emptyset \to 1, v \colon 2 \to \mathbf{2}, w \colon \mathbf{2}' \to \mathbf{2}\}$, where $2$ is the discrete category with two objects, $\mathbf{2}$ is the category with two objects $0, 1$ and a single non-identity morphism $0 \to 1$, and $\mathbf{2}'$ is the two-object category with two objects $0, 1$ and two non-identity morphisms $0 \to 1$.\par
We remark that $\mathscr{F}$-liftings offer a lot of flexibility in specification: For example if we want to get rid of the uniqueness assumption on the 2-cell $\varphi$ we can simply omit $w$ from $\mathcal{J}$.
\item\label{ex:liftingsGrpd} If $\mathcal{V} = \mathbf{Grpd}$, $\mathscr{F} = (\text{injective on objects}, \text{surjective equivalences})$, then the $\mathscr{F}$-lifting property is almost the same as in Example (\ref{ex:liftingsCat}) with the only difference being that all the 2-cells are invertible. Note that the weak factorization system $\mathscr{F}$ is cofibrantly generated by the set ${\mathcal{J} = \{u \colon \emptyset \to 1, v \colon 2 \to \mathbf{2}_g, w \colon \mathbf{2}'_g \to \mathbf{2}_g\}}$ that is almost the same as in Example (\ref{ex:liftingsCat}) with the only difference being that to each non-invertible morphism in $\mathbf{2}$ and $\mathbf{2}'$ we freely add its inverse, and in this way we obtain the groupoidal reflections $\mathbf{2}_g$ and $\mathbf{2}'_g$.
\item\label{ex:liftingsCh} If $\mathcal{V} = \mathbf{Ch}$, then we can choose $\mathcal{J} = \{S^{n - 1} \hookrightarrow D^n \mid n \in \mathbb{Z}\}$, $\mathscr{F} = ({}^\pitchfork(\mathcal{J}^\pitchfork), \mathcal{J}^\pitchfork)$, in which case $f \overset{\mathscr{F}}{\mathrel{\pitchfork}} k$ if and only if $\langle f,k \rangle$ is a surjective quasi-isomorphism. See \mbox{\cite[Proposition 2.3.4, Proposition 2.3.5]{Hov99}}.
\item\label{ex:liftingsSSet} Suppose that $\mathcal{V} = \mathbf{SSet}$ and let $\mathcal{J} = \{\partial \Delta^n \hookrightarrow \Delta^n \mid n \geq 0\}$, $\mathscr{F} = ({}^\pitchfork(\mathcal{J}^\pitchfork), \mathcal{J}^\pitchfork)$. Then $f \overset{\mathscr{F}}{\mathrel{\pitchfork}} k$ iff $\langle f,k \rangle$ is both a Kan fibration and a weak homotopy equivalence.
\item\label{ex:liftingsVIso} In this example we will show that the enriched lifting property in \cite{LW14} can be obtained as an $\mathscr{F}$-lifting property. Consider $\mathscr{F} = (\text{all maps}, \text{isomorphisms})$, which is obviously a weak factorization system on $\mathcal{V}_0$. Then $f \overset{\mathscr{F}}{\mathrel{\pitchfork}} k$ iff $\langle f,k \rangle$ is an isomorphism.\par
Now we will show that under some assumptions on $\mathcal{V}_0$ we can cofibrantly generate the weak factorization system $\mathscr{F}$. Assume that $\lambda$ is a regular cardinal and $\mathcal{V}_0$ is a locally $\lambda$-presentable category with a set $\mathcal{V}_\lambda$ of $\lambda$-presentable objects that form a strong generator. Consider $$\mathcal{J} = \{u_V \colon \emptyset \to V \mid V \in \mathcal{V}_\lambda\} \cup \{\nabla_V \colon V + V \to V \mid V \in \mathcal{V}_\lambda\},$$ where $\emptyset$ is the initial object and $\nabla_V = (\mathrm{id}_V,\mathrm{id}_V)$. Since $\mathcal{V}_\lambda$ is a strong generator, a morphism $v$ in $\mathcal{V}_0$ is an isomorphism iff $\mathcal{V}_0(V,v)$ is an isomorphism in $\mathbf{Set}$ for each $V \in \mathcal{V}_\lambda$. Now it suffices to notice that $u_V \mathrel{\pitchfork} v$ iff $\mathcal{V}_0(V,v)$ is surjective, and $\nabla_V \mathrel{\pitchfork} v$ iff $\mathcal{V}_0(V,v)$ is injective.
\item Suppose that $\mathscr{R}$ is the class of split epimorphisms in $\mathcal{V}_0$ and $\mathscr{L}$ is the class of retracts of binary coproduct injections in $\mathcal{V}_0$. By \cite[Proposition 2.6]{RT07}, $\mathscr{F} = (\mathscr{L}, \mathscr{R})$ is a weak factorization system. Then the $\mathscr{F}$-lifting property coincides with the $\mathcal{V}$\nobreakdash-enriched lifting property in the sense of \cite[Definition 13.3.1]{Rie14}. Finally, we remark that all the previous examples were cofibrantly generated, whereas this example is not necessarily cofibrantly generated.
\end{enumerate}
\endexamples
\definition
An \emph{enriched weak $\mathscr{F}$-factorization system} on $\mathcal{K}$ is a pair $(\mathcal{L}, \mathcal{R})$ of classes of morphisms in $\mathcal{K}_0$ such that each morphism $h \colon A \to B$ in $\mathcal{K}_0$ has a factorization $h = g \cdot f$ such that $f \in \mathcal{L}$, $g \in \mathcal{R}$, and furthermore $\mathcal{L} = {}^{\overset{\mathscr{F}}{\pitchfork}}\mathcal{R}$, $\mathcal{R} = \mathcal{L}^{\overset{\mathscr{F}}{\pitchfork}}$.
\enddefinition
\definition
We say that a class $\mathscr{S}$ of morphisms in $\mathcal{V}_0$ is \emph{stable under corners} if whenever $u, v \in \mathscr{S}$, then $u \mathrel{\Box} v \in \mathscr{S}$. (Recall Remark \ref{rem:Vnab}.)
\enddefinition
In the remainder of this section we will assume that $\mathcal{K}$ is a copowered $\mathcal{V}$-category admitting pushouts of the form $\mathrm{dom}(u \mathrel{\Box} f)$ from Definition \ref{def:uodotf}.
\rem\label{rem:enabadjlift}
The following equivalence holds: $$\text{$u \mathrel{\pitchfork} \langle f,k \rangle$ if and only if $u \mathrel{\Box} f \mathrel{\pitchfork} k$},$$ where $u$ is a morphism in $\mathcal{V}_0$ and $f$, $k$ are morphisms in $\mathcal{K}_0$. Indeed, this follows from Remark \ref{rem:adj} and Remark \ref{rem:enabadj}.
\endrem
\rem\label{rem:compatible}
From Remark \ref{rem:enabadjlift} we immediately conclude that for all morphisms $f$, $k$ in $\mathcal{K}_0$ the following equivalence holds:
$$\text{$u \mathrel{\Box} f \mathrel{\pitchfork} k$ holds for all $u \in \mathscr{L}$ if and only if $f \overset{\mathscr{F}}{\mathrel{\pitchfork}} k$}.$$
\endrem
\lemma\label{lem:strnabla}
Suppose that $\mathcal{J}$ is a class of morphisms in $\mathcal{V}_0$, and that the following implication holds: \begin{equation*}\text{If $u_1, u_2 \in \mathcal{J}$, then $u_1 \mathrel{\Box} u_2$ is in ${}^{\pitchfork}(\mathcal{J}^{\pitchfork})$.}\end{equation*} Then ${}^{\pitchfork}(\mathcal{J}^{\pitchfork})$ is stable under corners.
\endlemma
\pf
Suppose that $u_1 \in {}^{\pitchfork}(\mathcal{J}^{\pitchfork})$, $u_2 \in \mathcal{J}$, $v \in \mathcal{J}^\pitchfork$. By assumption we know that for each $u \in \mathcal{J}$: $u \mathrel{\Box} u_2 \mathrel{\pitchfork} v$, which is equivalent to $u \mathrel{\pitchfork} \langle u_2, v \rangle$ by Remark \ref{rem:enabadjlift}. Since this holds for all $u \in \mathcal{J}$, we get that $u_1 \mathrel{\pitchfork} \langle u_2, v \rangle$, and thus $u_1 \mathrel{\Box} u_2 \mathrel{\pitchfork} v$.\par
Now suppose that $u_1 \in {}^{\pitchfork}(\mathcal{J}^{\pitchfork})$, $u_2 \in {}^{\pitchfork}(\mathcal{J}^{\pitchfork})$, $v \in \mathcal{J}^\pitchfork$. By the previous paragraph, we know that for each $u \in \mathcal{J}$: $u_1 \mathrel{\Box} u \mathrel{\pitchfork} v$, which is equivalent to $u \mathrel{\pitchfork} \langle u_1, v \rangle$ by Remark \ref{rem:enabadjlift}. Since this holds for all $u \in \mathcal{J}$, we get that $u_2 \mathrel{\pitchfork} \langle u_1, v \rangle$, and thus $u_1 \mathrel{\Box} u_2 \mathrel{\pitchfork} v$. Note that we used the symmetry of the monoidal structure on $\mathcal{V}$.
\epf
\rem\label{rem:unitlaw}
Stability under corners forms part of the definition of a monoidal model category \cite[Definition 4.2.6]{Hov99}, and thus cofibrations and trivial cofibrations in a monoidal model category are stable under corners. The definition of a monoidal model category also involves a unit condition, and for our purposes of the enriched small object argument we do not need to assume any kind of a unit condition. However, the following interesting proposition characterizing self-enrichment of weak factorization systems in terms of a unit condition was suggested to the author by the anonymous referee.
\endrem
\proposition\label{prop:selfenrichment}
Suppose that $\mathscr{L}$ is stable under corners. Then $\mathscr{F}$ is an enriched weak $\mathscr{F}$-factorization system on $\mathcal{V}$ if and only if $!_I \colon \emptyset \to I$ belongs to $\mathscr{L}$.
\endproposition
\pf
Suppose that $!_I \in \mathscr{L}$. We know that $\mathscr{L}^\pitchfork = \mathscr{R}$, $\mathscr{L} = {}^{\pitchfork}\mathscr{R}$, and we want to show that $\mathscr{L} = {}^{\overset{\mathscr{F}}{\pitchfork}}\mathscr{R}$, $\mathscr{R} = \mathscr{L}^{\overset{\mathscr{F}}{\pitchfork}}$. We will show the first equality, the verification of the second equality is analogous. Suppose that $v \in {}^{\overset{\mathscr{F}}{\pitchfork}}\mathscr{R}$, i.e.\ for all $w \in \mathscr{R}$, $u \in \mathscr{L}$: $u \mathrel{\Box} v \mathrel{\pitchfork} w$. Thus in particular $!_I \mathrel{\Box} v \mathrel{\pitchfork} w$. Now it suffices to notice that $!_I \mathrel{\Box} v \cong v$ and we conclude ${}^{\overset{\mathscr{F}}{\pitchfork}}\mathscr{R} \subseteq \mathscr{L}$. Suppose that $v \in \mathscr{L}$. We want to show that $u \mathrel{\Box} v \mathrel{\pitchfork} w$ for all $u \in \mathscr{L}$, $w \in \mathscr{R}$. This in fact follows immediately from stability of $\mathscr{L}$ under corners, and thus we conclude that $\mathscr{L} \subseteq {}^{\overset{\mathscr{F}}{\pitchfork}}\mathscr{R}$.\par
On the other hand, suppose that $\mathscr{F}$ is an enriched weak $\mathscr{F}$-factorization system on $\mathcal{V}$. We will show that $!_I \in \mathscr{L}$. We know that $u \mathrel{\pitchfork} w$ for all $u \in \mathscr{L}$, $w \in \mathscr{R}$. Since $u \mathrel{\Box}\ !_I \cong u$, we obtain $u \mathrel{\Box}\ !_I \mathrel{\pitchfork} w$ for all $u \in \mathscr{L}$, $w \in \mathscr{R}$. Using $\mathscr{L} = {}^{\overset{\mathscr{F}}{\pitchfork}}\mathscr{R}$, we conclude $!_I \in \mathscr{L}$.
\epf

\section{Stability properties}\label{sec:stabilityproperties}

Here we establish the stability properties that are required in order to conclude that our small object argument generates an enriched weak $\mathscr{F}$-factorization system. Throughout the section we will assume that $\mathcal{I}$ is a class of morphisms in $\mathcal{K}_0$. We will be considering pushouts and transfinite composites in $\mathcal{K}$ as a $\mathcal{V}$-category: This means colimits in $\mathcal{K}_0$ that are sent by each representable functor $\mathcal{K}(-,K) \colon \mathcal{K}_0^\mathrm{op} \to \mathcal{V}_0$ to a limit in $\mathcal{V}_0$.
\proposition\label{prop:pushoutcl}
The class ${}^{\overset{\mathscr{F}}{\pitchfork}}\mathcal{I}$ is stable under pushouts in $\mathcal{K}$.
\endproposition
\pf
Suppose that the following square is a pushout in $\mathcal{K}$.
\[\begin{tikzcd}
	A & K \\
	B & L
	\arrow["f"', from=1-1, to=2-1]
	\arrow["g", from=1-1, to=1-2]
	\arrow["{f'}", from=1-2, to=2-2]
	\arrow["{g'}"', from=2-1, to=2-2]
\end{tikzcd}\]
We will show that if $f \in {}^{\overset{\mathscr{F}}{\pitchfork}}\mathcal{I}$, then $f' \in {}^{\overset{\mathscr{F}}{\pitchfork}}\mathcal{I}$.\par
Let $k \colon C \to D$ be in $\mathcal{I}$. We have induced morphisms $\langle f,k \rangle$, $\langle f',k \rangle$ in the notation of Definition \ref{def:inducede}. By assumption, $\langle f,k \rangle \in \mathscr{R}$. We want to show that $\langle f', k \rangle \in \mathscr{R}$. Recalling Definition \ref{def:sq}, there exists a morphism $\mathrm{Sq}((g,g'),k) \colon \mathrm{Sq}(f',k) \to \mathrm{Sq}(f,k)$ such that $p_1 \cdot \mathrm{Sq}((g,g'),k) = \mathcal{K}(g,C) \cdot p'_1$ and $p_2 \cdot \mathrm{Sq}((g,g'),k) = \mathcal{K}(g',D) \cdot p'_2$.
Now we will show that the square
\[\begin{tikzcd}[column sep=3em]
	{\mathcal{K}(L,C)} & {\mathcal{K}(B,C)} \\
	{\mathrm{Sq}(f',k)} & {\mathrm{Sq}(f,k)}
	\arrow["{\langle f',k \rangle}"', from=1-1, to=2-1]
	\arrow["{\mathcal{K}(g',C)}", from=1-1, to=1-2]
	\arrow["{\langle f,k \rangle}", from=1-2, to=2-2]
	\arrow["{\mathrm{Sq}((g,g'),k)}"', from=2-1, to=2-2]
\end{tikzcd}\]
is a pullback and this will finish the proof because $\mathscr{R}$ is stable under pullbacks.\par
Consider the following rectangle.
\[\begin{tikzcd}
	{\mathrm{Sq}(f',k)} & {\mathcal{K}(L,D)} & {\mathcal{K}(B,D)} \\
	{\mathcal{K}(K,C)} & {\mathcal{K}(K,D)} & {\mathcal{K}(A,D)}
	\arrow["{p'_1}"', from=1-1, to=2-1]
	\arrow["{p'_2}", from=1-1, to=1-2]
	\arrow["{\mathcal{K}(K,k)}"', from=2-1, to=2-2]
	\arrow["{\mathcal{K}(f',D)}", from=1-2, to=2-2]
	\arrow["{\mathcal{K}(g,D)}"', from=2-2, to=2-3]
	\arrow["{\mathcal{K}(g',D)}", from=1-2, to=1-3]
	\arrow["{\mathcal{K}(f,D)}", from=1-3, to=2-3]
\end{tikzcd}\]
The left square is a pullback by definition of $\mathrm{Sq}(f',k)$. The right square is a pullback, since $\mathcal{K}(-,D) \colon \mathcal{K}_0^{\mathrm{op}} \to \mathcal{V}_0$ preserves limits in $\mathcal{K}^{\mathrm{op}}$. Thus, the rectangle is a pullback by the pasting law for pullbacks.\par
Now consider the following rectangle.
\[\begin{tikzcd}[column sep=3em]
	{\mathrm{Sq}(f',k)} & {\mathrm{Sq}(f,k)} & {\mathcal{K}(B,D)} \\
	{\mathcal{K}(K,C)} & {\mathcal{K}(A,C)} & {\mathcal{K}(A,D)}
	\arrow["{p'_1}"', from=1-1, to=2-1]
	\arrow["{\mathrm{Sq}((g,g'),k)}", from=1-1, to=1-2]
	\arrow["{\mathcal{K}(g,C)}"', from=2-1, to=2-2]
	\arrow["{p_1}", from=1-2, to=2-2]
	\arrow["{\mathcal{K}(A,k)}"', from=2-2, to=2-3]
	\arrow["{p_2}", from=1-2, to=1-3]
	\arrow["{\mathcal{K}(f,D)}", from=1-3, to=2-3]
\end{tikzcd}\]
The right square is a pullback by definition of $\mathrm{Sq}(f,k)$, and the rectangle is a pullback, since $p_2 \cdot \mathrm{Sq}((g,g'),k) = \mathcal{K}(g',D) \cdot p'_2$ and $\mathcal{K}(A,k) \cdot \mathcal{K}(g,C) = \mathcal{K}(g,D) \cdot \mathcal{K}(K,k)$. Hence, by using the pasting law, we get that the left square in the rectangle is a pullback too.\par
Finally, consider the following rectangle.
\[\begin{tikzcd}
	{\mathcal{K}(L,C)} & {\mathrm{Sq}(f',k)} & {\mathcal{K}(K,C)} \\
	{\mathcal{K}(B,C)} & {\mathrm{Sq}(f,k)} & {\mathcal{K}(A,C)}
	\arrow["{\mathcal{K}(g',C)}"', from=1-1, to=2-1]
	\arrow["{\langle f',k \rangle}", from=1-1, to=1-2]
	\arrow["{\langle f,k \rangle}"', from=2-1, to=2-2]
	\arrow["{\mathrm{Sq}((g,g'),k)}", from=1-2, to=2-2]
	\arrow["{p_1}"', from=2-2, to=2-3]
	\arrow["{p'_1}", from=1-2, to=1-3]
	\arrow["{\mathcal{K}(g,C)}", from=1-3, to=2-3]
\end{tikzcd}\]
The right square was shown above to be a pullback. Moreover, the rectangle is a pullback because by composing the horizontal sides we get the square
\[\begin{tikzcd}
	{\mathcal{K}(L,C)} & {\mathcal{K}(K,C)} \\
	{\mathcal{K}(B,C)} & {\mathcal{K}(A,C)}
	\arrow["{\mathcal{K}(g',C)}"', from=1-1, to=2-1]
	\arrow["{\mathcal{K}(g,C)}", from=1-2, to=2-2]
	\arrow["{\mathcal{K}(f',C)}", from=1-1, to=1-2]
	\arrow["{\mathcal{K}(f,C)}"', from=2-1, to=2-2]
\end{tikzcd}\]
and this square is a pullback, since $\mathcal{K}(-,C) \colon \mathcal{K}_0^{\mathrm{op}} \to \mathcal{V}_0$ preserves limits in $\mathcal{K}^{\mathrm{op}}$. Therefore, from the pasting law, we conclude that the left square in the rectangle is a pullback.
\epf
\proposition
The class ${}^{\overset{\mathscr{F}}{\pitchfork}}\mathcal{I}$ is stable under transfinite compositions in $\mathcal{K}$.
\endproposition
\pf
Suppose that $\alpha > 0$ is an ordinal. We will show that if $f_{\beta,\beta + 1} \colon A_\beta \to A_{\beta + 1}$, $\beta < \alpha$, are morphisms coming from a functor $A_{-} \colon \alpha \to \mathcal{K}_0$ such that each $f_{\beta, \beta + 1}$ belongs to ${}^{\overset{\mathscr{F}}{\pitchfork}}\mathcal{I}$, and for each limit ordinal $\gamma < \alpha$ the induced morphism $i_\gamma \colon \operatorname{colim}_{\delta < \gamma} A_\delta \to A_\gamma$ is an isomorphism, then their transfinite composition $f_{0,\alpha} \colon A_0 \to A_\alpha := \operatorname{colim}_{\beta < \alpha} A_\beta$ also belongs to ${}^{\overset{\mathscr{F}}{\pitchfork}}\mathcal{I}$.\par
The proof is by transfinite induction. \par
Base Case: If $\alpha = 1$, then the statement clearly holds: $f_{0,1} \in {}^{\overset{\mathscr{F}}{\pitchfork}}\mathcal{I}$ implies that $f_{0,1} \in {}^{\overset{\mathscr{F}}{\pitchfork}}\mathcal{I}$. \par
Successor Step: Suppose that $\alpha = \epsilon + 1 > 1$ and that the result holds for all non-zero ordinals less than $\alpha$. By inductive hypothesis we know that $f_{0,\epsilon} \colon A_0 \to A_\epsilon$ belongs to ${}^{\overset{\mathscr{F}}{\pitchfork}}\mathcal{I}$. Let $k \colon C \to D$ be in $\mathcal{I}$. To simplify notation, denote $A := A_0$, $A' := A_{\epsilon}$, $B := A_{\alpha}$, $f := f_{0,\epsilon}$, $f' := f_{\epsilon, \epsilon + 1}$. Note that $f_{0, \alpha} = f' \cdot f$. We have induced morphisms $\langle f,k \rangle, \langle f',k \rangle$, $\langle f' \cdot f,k \rangle$ from Definition \ref{def:inducede} whose respective pullback projections will be denoted $p_i$, $p'_i$, $p''_i$, where ${i \in \{1,2\}}$. We know that $\langle f,k \rangle \in \mathscr{R}$, $\langle f',k \rangle \in \mathscr{R}$, and we want to show that $\langle f' \cdot f,k \rangle \in \mathscr{R}$. Recalling Definition \ref{def:sq}, we have a morphism ${\mathrm{Sq}((f,\mathrm{id}_B),k) \colon \mathrm{Sq}(f',k) \to \mathrm{Sq}(f' \cdot f,k)}$ such that $p''_1 \cdot \mathrm{Sq}((f,\mathrm{id}_B),k) = \mathcal{K}(f,C) \cdot p'_1$ and $p''_2 \cdot \mathrm{Sq}((f,\mathrm{id}_B),k) = p'_2$, and also a morphism $\mathrm{Sq}((\mathrm{id}_A,f'),k) \colon \mathrm{Sq}(f' \cdot f,k)\to \mathrm{Sq}(f,k)$ such that $p_2 \cdot \mathrm{Sq}((\mathrm{id}_A,f'),k) = \mathcal{K}(f',D) \cdot p''_2$ and $p_1 \cdot \mathrm{Sq}((\mathrm{id}_A,f'),k) = p''_1$.
Now note that 
$$p''_2 \cdot \mathrm{Sq}((f,\mathrm{id}_B),k) \cdot \langle f',k \rangle = p'_2 \cdot \langle f',k \rangle = \mathcal{K}(B,k) = p''_2 \cdot \langle f' \cdot f,k \rangle$$
and
\begin{align*}
p''_1 \cdot \mathrm{Sq}((f,\mathrm{id}_B),k) \cdot \langle f',k \rangle &= \mathcal{K}(f,C) \cdot p'_1 \cdot \langle f',k \rangle \\ &= \mathcal{K}(f,C) \cdot \mathcal{K}(f',C) \\ &= \mathcal{K}(f' \cdot f,C) \\ &= p''_1 \cdot \langle f' \cdot f,k \rangle.
\end{align*}
Therefore $\langle f' \cdot f,k \rangle = \mathrm{Sq}((f,\mathrm{id}_B),k) \cdot \langle f',k \rangle$, since $(p''_1, p''_2)$ is a pullback. We will show that $\mathrm{Sq}((f,\mathrm{id}_B),k) \in \mathscr{R}$ and this will finish the proof of the successor step, since $\mathscr{R}$ is stable under compositions. In order to do that, we will show that the square
\[\begin{tikzcd}[column sep=3.6em]
	{\mathrm{Sq}(f',k)} & {\mathcal{K}(A',C)} \\
	{\mathrm{Sq}(f' \cdot f,k)} & {\mathrm{Sq}(f,k)}
	\arrow["{p'_1}", from=1-1, to=1-2]
	\arrow["{\mathrm{Sq}((f,\mathrm{id}_B),k)}"', from=1-1, to=2-1]
	\arrow["{\langle f,k \rangle}", from=1-2, to=2-2]
	\arrow["{\mathrm{Sq}((\mathrm{id}_A,f'),k) }"', from=2-1, to=2-2]
\end{tikzcd}\]
is a pullback and this will imply $\mathrm{Sq}((f,\mathrm{id}_B),k) \in \mathscr{R}$, since $\mathscr{R}$ is stable under pullbacks. Indeed, the following diagram commutes.
\[\begin{tikzcd}[column sep = 3.6em]
	{\mathrm{Sq}(f',k)} & {\mathcal{K}(A',C)} \\
	{\mathrm{Sq}(f' \cdot f,k)} & {\mathrm{Sq}(f,k)} & {\mathcal{K}(A,C)} \\
	{\mathcal{K}(B,D)} & {\mathcal{K}(A',D)} & {\mathcal{K}(A,D)}
	\arrow["{p'_1}", from=1-1, to=1-2]
	\arrow["{\mathrm{Sq}((f,\mathrm{id}_B),k)}"', from=1-1, to=2-1]
	\arrow["{\langle f,k \rangle}", from=1-2, to=2-2]
	\arrow["{\mathrm{Sq}((\mathrm{id}_A,f'),k) }", from=2-1, to=2-2]
	\arrow["{p''_2}"', from=2-1, to=3-1]
	\arrow["{p_1}", from=2-2, to=2-3]
	\arrow["{p_2}", from=2-2, to=3-2]
	\arrow["{\mathcal{K}(A,k)}", from=2-3, to=3-3]
	\arrow["{\mathcal{K}(f',D)}"', from=3-1, to=3-2]
	\arrow["{\mathcal{K}(f,D)}"', from=3-2, to=3-3]
\end{tikzcd}\]
Moreover, the horizontal rectangle is a pullback and the right square in the horizontal rectangle is a pullback, thus, by the pasting law, the left square in the horizontal rectangle is a pullback. Finally, the vertical rectangle is a pullback and the bottom square in the vertical rectangle is a pullback, hence, by the pasting law, the top square in the vertical rectangle is a pullback.
\par
Limit Step: Suppose that $\alpha > 0$ is a limit ordinal such that the result holds for all non-zero ordinals less than $\alpha$. Let $k \colon C \to D$ be in $\mathcal{I}$. We know that for each $\beta < \alpha$ the induced morphism $\langle f_{\beta, \beta + 1}, k \rangle$ from the following diagram belongs to $\mathscr{R}$.
\[\begin{tikzcd}[column sep=3em]
	{\mathcal{K}(A_{\beta + 1},C)} \\
	& {\mathrm{Sq}(f_{\beta,\beta + 1},k)} & {\mathcal{K}(A_\beta,C)} \\
	& {\mathcal{K}(A_{\beta + 1},D)} & {\mathcal{K}(A_\beta,D)}
	\arrow["{p_1^{\beta,\beta+1}}"', from=2-2, to=2-3]
	\arrow["{p_2^{\beta,\beta+1}}", from=2-2, to=3-2]
	\arrow["{\mathcal{K}(A_\beta,k)}", from=2-3, to=3-3]
	\arrow["{\mathcal{K}(f_{\beta,\beta+1},D)}"', from=3-2, to=3-3]
	\arrow["{\mathcal{K}(A_{\beta + 1},k)}"', from=1-1, to=3-2]
	\arrow["{\mathcal{K}(f_{\beta,\beta + 1},C)}", from=1-1, to=2-3]
	\arrow["{\langle f_{\beta,\beta + 1},k \rangle}"{description}, dashed, from=1-1, to=2-2]
\end{tikzcd}\] 
Let $(f_{\beta,\alpha} \colon A_\beta \to A_\alpha)_{\beta < \alpha}$ be a colimit of the diagram $\{f_{\beta, \beta + 1} \colon A_{\beta} \to A_{\beta + 1} \mid \beta < \alpha\}$. For simplicity, we will write the diagram as follows.
\[\begin{tikzcd}
	{A_0} & {A_1} & {A_2} & {A_3} & \cdots
	\arrow["{f_{0,1}}", from=1-1, to=1-2]
	\arrow["{f_{1,2}}", from=1-2, to=1-3]
	\arrow["{f_{2,3}}", from=1-3, to=1-4]
	\arrow["{f_{3,4}}", from=1-4, to=1-5]
\end{tikzcd}\]
Then $((f_{\beta,\alpha}, \mathrm{id}_{A_\alpha}) \colon f_{\beta,\alpha} \to \mathrm{id}_{A_\alpha})_{\beta < \alpha}$ is a colimit of the following diagram in $[\mathbf{2},\mathcal{K}]$, since it is a colimit in each component.
\[\begin{tikzcd}
	{A_0} & {A_1} & {A_2} & {A_3} & \cdots \\
	{A_\alpha} & {A_\alpha} & {A_\alpha} & {A_\alpha} & \dots
	\arrow["{f_{0,1}}", from=1-1, to=1-2]
	\arrow["{f_{1,2}}", from=1-2, to=1-3]
	\arrow["{f_{2,3}}", from=1-3, to=1-4]
	\arrow["{f_{3,4}}", from=1-4, to=1-5]
	\arrow["{f_{0,\alpha}}"', from=1-1, to=2-1]
	\arrow["{f_{1,\alpha}}"', from=1-2, to=2-2]
	\arrow["{\mathrm{id}_{A_\alpha}}"', from=2-1, to=2-2]
	\arrow["{f_{2,\alpha}}"', from=1-3, to=2-3]
	\arrow["{f_{3,\alpha}}"', from=1-4, to=2-4]
	\arrow["{\mathrm{id}_{A_\alpha}}"', from=2-2, to=2-3]
	\arrow["{\mathrm{id}_{A_\alpha}}"', from=2-3, to=2-4]
	\arrow["{\mathrm{id}_{A_\alpha}}"', from=2-4, to=2-5]
\end{tikzcd}\]
Recall from Definition \ref{def:sq} that $[\mathbf{2},\mathcal{K}]$ is a $\mathcal{V}$-category with $\mathrm{Sq}(g,h)$ being the hom-object for each pair $g, h \in [\mathbf{2},\mathcal{K}]$. Thus, $$\big(\mathrm{Sq}((f_{\beta,\alpha}, \mathrm{id}_{A_\alpha}),k) \colon \mathrm{Sq}(\mathrm{id}_{A_\alpha},k) \to \mathrm{Sq}(f_{\beta,\alpha},k)\big)_{\beta < \alpha}$$ is a limit of the diagram
\[\begin{tikzcd}[column sep=2.8em]
	{\mathrm{Sq}(f_{0,\alpha},k)} && {\mathrm{Sq}(f_{1,\alpha},k)} && {\mathrm{Sq}(f_{2,\alpha},k)} && \cdots
	\arrow["{\mathrm{Sq}((f_{0,1},\mathrm{id}_{A_\alpha}),k)}"', from=1-3, to=1-1]
	\arrow["{\mathrm{Sq}((f_{1,2},\mathrm{id}_{A_\alpha}),k)}"', from=1-5, to=1-3]
	\arrow["{\mathrm{Sq}((f_{2,3},\mathrm{id}_{A_\alpha}),k)}"', from=1-7, to=1-5]
\end{tikzcd}\]
in $\mathcal{V}_0$ because $\mathrm{Sq}(-,k) \colon [\mathbf{2},\mathcal{K}]_0^{\mathrm{op}} \to \mathcal{V}_0$ preserves limits in $[\mathbf{2},\mathcal{K}]^{\mathrm{op}}$. Also, note that there exists an isomorphism $\iota \colon \mathcal{K}(A_\alpha,C) \to \mathrm{Sq}(\mathrm{id}_{A_\alpha},k)$ that makes the following diagram commute because the second component is uniquely determined by the first component via composition with $k$.
\[\begin{tikzcd}
	& {\mathcal{K}(A_\alpha,C)} \\
	{\mathcal{K}(A_\alpha,C)} & {\mathrm{Sq}(\mathrm{id}_{A_\alpha},k) } & {\mathcal{K}(A_{\alpha},D)}
	\arrow["{p_1}", from=2-2, to=2-1]
	\arrow["{p_2}"', from=2-2, to=2-3]
	\arrow["{\mathcal{K}(A_\alpha,k)}", from=1-2, to=2-3]
	\arrow["{\mathrm{id}_{\mathcal{K}(A_\alpha,C)}}"', from=1-2, to=2-1]
	\arrow["\iota", from=1-2, to=2-2]
\end{tikzcd}\]
Thus, the transfinite cocomposition $\mathrm{Sq}((f_{0,\alpha}, \mathrm{id}_{A_\alpha}),k) \colon \mathrm{Sq}(\mathrm{id}_{A_\alpha},k) \to \mathrm{Sq}(f_{0,\alpha},k)$ is isomorphic to a morphism $\mathrm{Sq}((f_{0,\alpha}, \mathrm{id}_{A_\alpha}),k) \cdot \iota \colon \mathcal{K}(A_\alpha,C) \to \mathrm{Sq}(f_{0,\alpha},k)$, and this morphism is the induced map $\langle f_{0,\alpha},k \rangle$ in the following diagram by the uniqueness of the induced map.
\[\begin{tikzcd}
	{\mathcal{K}(A_{\alpha},C)} \\
	& {\mathrm{Sq}(f_{0,\alpha},k)} & {\mathcal{K}(A_0,C)} \\
	& {\mathcal{K}(A_{\alpha},D)} & {\mathcal{K}(A_0,D)}
	\arrow["{p_1^{0,\alpha}}"', from=2-2, to=2-3]
	\arrow["{p_2^{0,\alpha}}", from=2-2, to=3-2]
	\arrow["{\mathcal{K}(A_0,k)}", from=2-3, to=3-3]
	\arrow["{\mathcal{K}(f_{0,\alpha},D)}"', from=3-2, to=3-3]
	\arrow["{\mathcal{K}(A_{\alpha},k)}"', from=1-1, to=3-2]
	\arrow["{\mathcal{K}(f_{0,\alpha},C)}", from=1-1, to=2-3]
	\arrow["{\langle f_{0,\alpha},k \rangle}"{description}, dashed, from=1-1, to=2-2]
\end{tikzcd}\]
Indeed, 
$$p_1^{0,\alpha} \cdot \mathrm{Sq}((f_{0,\alpha}, \mathrm{id}_{A_\alpha}),k) \cdot \iota = \mathcal{K}(f_{0,\alpha},C) \cdot p_1 \cdot \iota = \mathcal{K}(f_{0,\alpha},C) \cdot \mathrm{id}_{\mathcal{K}(A_\alpha,C)} = \mathcal{K}(f_{0,\alpha},C)$$
and
$$p_2^{0,\alpha} \cdot \mathrm{Sq}((f_{0,\alpha}, \mathrm{id}_{A_\alpha}),k) \cdot \iota = \mathrm{id}_{A_\alpha} \cdot p_2 \cdot \iota = p_2 \cdot \iota = \mathcal{K}(A_\alpha,k).$$
Therefore, $\mathrm{Sq}((f_{0,\alpha}, \mathrm{id}_{A_\alpha}),k) \cdot \iota = \langle f_{0,\alpha},k \rangle$. In order to finish the proof we need to show that $\langle f_{0,\alpha},k \rangle$ belongs to $\mathscr{R}$. To achieve that, recall that $\mathscr{R}$ is stable under transfinite cocompositions and isomorphisms, and thus it suffices to show that for each $\beta < \alpha$ the morphism $\mathrm{Sq}((f_{\beta,\beta+1},\mathrm{id}_{A_\alpha}),k)$ belongs to $\mathscr{R}$. Let $\beta < \alpha$. We will show that the square
\[\begin{tikzcd}[column sep=5.3em]
	{\mathrm{Sq}(f_{\beta+1,\alpha},k)} & {\mathcal{K}(A_{\beta + 1},C)} \\
	{\mathrm{Sq}(f_{\beta,\alpha},k)} & {\mathrm{Sq}(f_{\beta,\beta,+1},k)}
	\arrow["{p_1^{\beta+1,\alpha}}", from=1-1, to=1-2]
	\arrow["{\mathrm{Sq}((f_{\beta,\beta+1},\mathrm{id}_{A_\alpha}),k)}"', from=1-1, to=2-1]
	\arrow["{\langle f_{\beta,\beta+1},k \rangle}", from=1-2, to=2-2]
	\arrow["{\mathrm{Sq}((\mathrm{id}_{A_\beta},f_{\beta+1,\alpha}),k)}"', from=2-1, to=2-2]
\end{tikzcd}\]
is a pullback and this will imply that $\mathrm{Sq}((f_{\beta,\beta+1},\mathrm{id}_{A_\alpha}),k)$ belongs to $\mathscr{R}$ because $\mathscr{R}$ is stable under pullbacks. Indeed, the following diagram commutes.
\[\begin{tikzcd}[column sep = 3.5em]
	{\mathrm{Sq}(f_{\beta +1,\alpha},k)} & {\mathcal{K}(A_{\beta + 1},C)} \\
	{\mathrm{Sq}(f_{\beta,\alpha},k)} & {\mathrm{Sq}(f_{\beta,\beta+1},k)} & {\mathcal{K}(A_\beta,C)} \\
	{\mathcal{K}(A_{\alpha},D)} & {\mathcal{K}(A_{\beta + 1},D)} & {\mathcal{K}(A_\beta,D)}
	\arrow["{p_1^{\beta + 1,\alpha}}", from=1-1, to=1-2]
	\arrow["{\mathrm{Sq}((f_{\beta,\beta+1},\mathrm{id}_{A_\alpha}),k)}"', from=1-1, to=2-1]
	\arrow["{\langle f_{\beta,\beta + 1},k \rangle}", from=1-2, to=2-2]
	\arrow["{\mathrm{Sq}((\mathrm{id}_{A_\beta},f_{\beta+1,\alpha}),k)}", from=2-1, to=2-2]
	\arrow["{p_2^{\beta,\alpha}}"', from=2-1, to=3-1]
	\arrow["{p_1^{\beta,\beta+1}}", from=2-2, to=2-3]
	\arrow["{p_2^{\beta,\beta+1}}", from=2-2, to=3-2]
	\arrow["{\mathcal{K}(A_\beta,k)}", from=2-3, to=3-3]
	\arrow["{\mathcal{K}(f_{\beta + 1,\alpha},D)}"', from=3-1, to=3-2]
	\arrow["{\mathcal{K}(f_{\beta, \beta + 1},D)}"', from=3-2, to=3-3]
\end{tikzcd}\]
Moreover, the horizontal rectangle is a pullback and the right square in the horizontal rectangle is a pullback, thus, by the pasting law, the left square in the horizontal rectangle is a pullback. Finally, the vertical rectangle is a pullback and the bottom square in the vertical rectangle is a pullback, hence, by the pasting law, the top square in the vertical rectangle is a pullback.
\epf
In the remainder of this section we will assume that $\mathcal{K}$ is a copowered $\mathcal{V}$-category admitting pushouts of the form $\mathrm{dom}(u \mathrel{\Box} f)$ from Definition \ref{def:uodotf}.
\proposition\label{prop:copowpush}
Suppose that $\mathscr{L}$ is stable under corners. Moreover, suppose that we have $u \colon U \to V$ in $\mathscr{L}$ and $f \colon A \to B$ in $\mathcal{K}_0$ such that $f \in {}^{\overset{\mathscr{F}}{\pitchfork}}\mathcal{I}$. Then $u \mathrel{\Box} f \in {}^{\overset{\mathscr{F}}{\pitchfork}}\mathcal{I}$.
\endproposition
\pf
In order to show that $u \mathrel{\Box} f \in {}^{\overset{\mathscr{F}}{\pitchfork}}\mathcal{I}$ it suffices to show that ${v \mathrel{\Box} (u \mathrel{\Box} f) \in {}^{\pitchfork}\mathcal{I}}$ for each $v \in \mathscr{L}$, see Remark \ref{rem:compatible}. From Equation (\ref{eqn:nabassoc}) we know that $v \mathrel{\Box} (u \mathrel{\Box} f)$ is isomorphic to $(v \mathrel{\Box} u) \mathrel{\Box} f$. Note that $v \mathrel{\Box} u$ is in $\mathscr{L}$, since $\mathscr{L}$ is stable under corners. Furthermore, $f$ is in ${}^{\overset{\mathscr{F}}{\pitchfork}}\mathcal{I}$. Hence, recalling Remark \ref{rem:compatible}, we obtain ${(v \mathrel{\Box} u) \mathrel{\Box} f \in {}^{\pitchfork}\mathcal{I}}$, and thus $v \mathrel{\Box} (u \mathrel{\Box} f) \in {}^{\pitchfork}\mathcal{I}$.
\epf
\corollary\label{cor:copowpush}
Suppose that $\mathscr{L}$ is stable under corners. Moreover, suppose that we have a span
\[\begin{tikzcd}
	{\mathrm{dom}(u \mathrel{\Box} f)} & K \\
	{V \odot B}
	\arrow["{u \mathrel{\Box} f}"', from=1-1, to=2-1]
	\arrow["{(h,g)}", from=1-1, to=1-2]
\end{tikzcd}\]
in $\mathcal{K}_0$, where $u \colon U \to V$ is in $\mathscr{L}$, $f \colon A \to B$, $g \colon V \odot A \to K$, and $h \colon U \odot B \to K$ are morphisms in $\mathcal{K}_0$, and that the following square is a pushout in $\mathcal{K}$.
\begin{equation}\label{cmd:cowpushfirstdef}
\begin{tikzcd}
	{\mathrm{dom}(u \mathrel{\Box} f)} & K \\
	{V \odot B} & L
	\arrow["{u \mathrel{\Box} f}"', from=1-1, to=2-1]
	\arrow["{(h,g)}", from=1-1, to=1-2]
	\arrow["{f'}", from=1-2, to=2-2]
	\arrow["{g'}"', from=2-1, to=2-2]
\end{tikzcd}
\end{equation}
If $f \in {}^{\overset{\mathscr{F}}{\pitchfork}}\mathcal{I}$, then $f' \in {}^{\overset{\mathscr{F}}{\pitchfork}}\mathcal{I}$.
\endcorollary
\pf
From Proposition \ref{prop:copowpush} we get that $u \mathrel{\Box} f \in {}^{\overset{\mathscr{F}}{\pitchfork}}\mathcal{I}$. Now it suffices to recall that ${}^{\overset{\mathscr{F}}{\pitchfork}}\mathcal{I}$ is stable under pushouts, see Proposition \ref{prop:pushoutcl}.
\epf
\definition
We call the pushout (\ref{cmd:cowpushfirstdef}) a \emph{copowered pushout of $f, g, h$ relative to $u$}.\par
Furthermore, recall that that $\mathcal{I}$ is a class of morphisms in $\mathcal{K}_0$ and suppose that $\mathcal{J}$ is a class of morphisms in $\mathcal{V}_0$. Then we say that a morphism $f'$ in $\mathcal{K}_0$ is \emph{a copowered pushout of a morphism from $\mathcal{I}$ relative to $\mathcal{J}$} if there exist morphisms $u \in \mathcal{J}$, $f \in \mathcal{I}$ and morphisms $g, h, g'$ in $\mathcal{K}_0$ such that the square (\ref{cmd:cowpushfirstdef}) is a pushout in $\mathcal{K}$.
\enddefinition
\rem
Note that all of the information concerning the span from Corollary \ref{cor:copowpush} is encoded in the following diagram
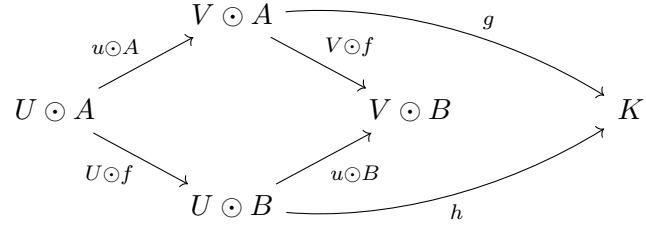
\[\begin{tikzcd}
	& {V \odot A} \\
	{U \odot A} && {V \odot B} && K \\
	& {U \odot B}
	\arrow["h"', curve={height=18pt}, from=3-2, to=2-5]
	\arrow["{U \odot f}"', from=2-1, to=3-2]
	\arrow["{V \odot f }", from=1-2, to=2-3]
	\arrow["g"{pos=0.6}, curve={height=-18pt}, from=1-2, to=2-5]
	\arrow["{u \odot A}", from=2-1, to=1-2]
	\arrow["{u \odot B}"', from=3-2, to=2-3]
\end{tikzcd}\]
whose colimit is the copowered pushout as displayed below.
\[\begin{tikzcd}[sep=2.8em]
	& {V \odot A} \\
	{U \odot A} && {V \odot B} & L & K \\
	& {U \odot B}
	\arrow["h"', curve={height=18pt}, from=3-2, to=2-5]
	\arrow["{U \odot f}"', from=2-1, to=3-2]
	\arrow["{V \odot f }", from=1-2, to=2-3]
	\arrow["g"{pos=0.6}, curve={height=-18pt}, from=1-2, to=2-5]
	\arrow["{f'}"', from=2-5, to=2-4]
	\arrow["{g'}"{pos=0.6}, from=2-3, to=2-4]
	\arrow["{u \odot A}", from=2-1, to=1-2]
	\arrow["{u \odot B}"', from=3-2, to=2-3]
\end{tikzcd}\]
The rephrasing, together with the bijective correspondence from Remark \ref{rem:copbij}, is useful for understanding Examples \ref{ex:cpwrpush}.
\endrem
\examples\label{ex:cpwrpush} In the following examples we will describe the part of the universal property of a copowered pushout that takes place in $\mathcal{K}_0$. In each example we assume that $\mathscr{F} = ({}^\pitchfork(\mathcal{J}^\pitchfork), \mathcal{J}^\pitchfork)$.
\begin{enumerate}[(1)]
\item In the case $\mathcal{V} = \mathbf{Set}$, $\mathcal{J} = \{u \colon \emptyset \to 1\}$, a copowered pushout of $f, g, h$ relative to $u$ is the same as a pushout of $f, g$. Note that $h$ is the unique morphism $\emptyset \to K$.\par
Moreover, in the case $\mathcal{V} = \mathbf{Set}$, $\mathcal{J} = \{u \colon \emptyset \to 1, v \colon 2 \to 1\}$, a copowered pushout of $f, g, h$ relative to $v$ takes as input the following diagram
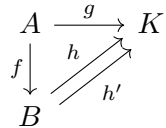
\[\begin{tikzcd}
	A & K \\
	B
	\arrow["g", from=1-1, to=1-2]
	\arrow["f"', from=1-1, to=2-1]
	\arrow["{h'}"', shift right, from=2-1, to=1-2]
	\arrow["h", shift left, from=2-1, to=1-2]
\end{tikzcd}\]
in which $h \cdot f = h' \cdot f = g$ and returns a morphism $c \colon K \to L$ such that $c \cdot h = c \cdot h'$ with the following universal property: For each morphism $d \colon K \to D$ satisfying $d \cdot h = d \cdot h'$ there exists a unique morphism $p \colon L \to D$ such that $p \cdot c = d$. Note that $c$ is the coequalizer of $h$ and $h'$.
\item\label{ex:cpwrpushCat} Suppose that $\mathcal{V} = \mathbf{Cat}$ and $\mathcal{J} = \{u \colon \emptyset \to 1, v \colon 2 \to \mathbf{2}, w \colon \mathbf{2}' \to \mathbf{2}\}$.
\begin{enumerate}
\item A copowered pushout of $f, g, h$ relative to $u$ is the same as a pushout of $f, g$.
\item A copowered pushout of $f, g, h$ relative to $v$ takes as input cells in $\mathcal{K}$ depicted in the diagram below
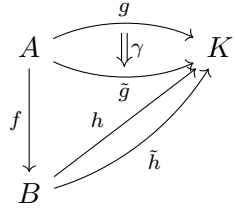
\[\begin{tikzcd}
	A && K \\
	\\
	B
	\arrow[""{name=0, anchor=center, inner sep=0}, "g", curve={height=-12pt}, from=1-1, to=1-3]
	\arrow[""{name=1, anchor=center, inner sep=0}, "{\tilde{g}}"', curve={height=12pt}, from=1-1, to=1-3]
	\arrow["f"', from=1-1, to=3-1]
	\arrow["h"{pos=0.4}, from=3-1, to=1-3]
	\arrow["{\tilde{h}}"', curve={height=12pt}, from=3-1, to=1-3]
	\arrow["\gamma", shorten <=3pt, shorten >=3pt, Rightarrow, from=0, to=1]
\end{tikzcd}\]
in which $h \cdot f = g$, $\tilde{h} \cdot f = \tilde{g}$, and returns a cocone $j, \tilde{j} \colon B \to L$, $\delta \colon j \Rightarrow \tilde{j}$, $i \colon K \to L$ such that $j \cdot f = i \cdot g$, $\tilde{j} \cdot f = i \cdot \tilde{g}$, $j = i \cdot h$, $\tilde{j} = i \cdot \tilde{h}$, and $\delta \ast f = i \ast \gamma$ with the following universal property: For each other such compatible cocone $r \colon K \to D$, $s, \tilde{s} \colon B \to D$, $\epsilon \colon s \Rightarrow \tilde{s}$ there exists a unique 1-cell $p \colon L \to D$ such that $p \ast \delta = \epsilon$ and $p \cdot i = r$.
\item A copowered pushout of $f, g, h$ relative to $w$ takes as input cells in $\mathcal{K}$ depicted in the diagram below
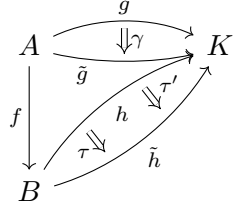
\[\begin{tikzcd}
	A && K \\
	\\
	B
	\arrow["f"', from=1-1, to=3-1]
	\arrow[""{name=0, anchor=center, inner sep=0}, "g", curve={height=-12pt}, from=1-1, to=1-3]
	\arrow[""{name=1, anchor=center, inner sep=0}, "h"', curve={height=-12pt}, from=3-1, to=1-3]
	\arrow[""{name=2, anchor=center, inner sep=0}, "{\tilde{h}}"', curve={height=12pt}, from=3-1, to=1-3]
	\arrow[""{name=3, anchor=center, inner sep=0}, "{\tilde{g}}"'{pos=0.3}, curve={height=6pt}, from=1-1, to=1-3]
	\arrow["\tau"', shift right=5, shorten <=4pt, shorten >=4pt, Rightarrow, from=1, to=2]
	\arrow["{\tau'}", shift left=5, shorten <=4pt, shorten >=4pt, Rightarrow, from=1, to=2]
	\arrow["\gamma", shorten <=2pt, shorten >=2pt, Rightarrow, from=0, to=3]
\end{tikzcd}\]
in which $h \cdot f = g$, $\tilde{h} \cdot f = \tilde{g}$, $\tau \ast f = \gamma$, $\tau' \ast f = \gamma$, and returns a 1-cell $c \colon K \to L$ such that $c \ast \tau = c \ast \tau'$ with the following universal property: For any 1-cell $d \colon K \to D$ satisfying $d \ast \tau = d \ast \tau'$ there exists a unique 1-cell $p \colon L \to D$ such that $p \cdot c = d$. Note that $c$ is the coequifier of $\tau$ and $\tau'$.
\end{enumerate}
\item\label{ex:cpwrpushGrpd} If $\mathcal{V} = \mathbf{Grpd}$, $\mathcal{J} = \{u \colon \emptyset \to 1, v \colon 2 \to \mathbf{2}_g, w \colon \mathbf{2}'_g \to \mathbf{2}_g\}$, then copowered pushouts are almost the same as in Example (\ref{ex:cpwrpushCat}) with the only difference being that all the 2-cells are invertible.
\item In the case $\mathcal{V} = \mathbf{Ch}$, $\mathcal{J} = \{u_n \colon S^{n - 1} \hookrightarrow D^n \mid n \in \mathbb{Z}\}$, a copowered pushout of $f, g, h$ relative to $u_n$ takes as input an $n$-chain $g \in \mathcal{K}(A,K)_n$, a morphism $f \colon A \to B$ in $\mathcal{K}_0$, and an $(n - 1)$-cycle $h \in \mathcal{K}(B,K)_{n - 1}$ such that $h \cdot f = \partial_n(g)$, and it returns a morphism $f' \colon K \to L$ in $\mathcal{K}_0$ and an $n$-chain $g' \in \mathcal{K}(B,L)_n$ such that $f' \cdot h = \partial_n(g')$ and $f' \cdot g = g' \cdot f \in \mathcal{K}(A,K)_n$ with the following universal property: For each other such compatible pair consisting of a morphism $r \colon K \to D$ in $\mathcal{K}_0$ and an $n$-chain $s \in \mathcal{K}(B,D)_n$ there exists a unique morphism $p \colon L \to D$ in $\mathcal{K}_0$ such that $p \cdot f' = r$ and $p \cdot g' = s$.
\item Suppose that $\mathcal{V} = \mathbf{SSet}$ and $\mathcal{J} = \{u_n \colon \partial \Delta^n \hookrightarrow \Delta^n \mid n \geq 0\}$. A copowered pushout of $f, g, h$ relative to $u_n$ takes as input an $n$-simplex $g \in \mathcal{K}(A,K)_n$, a morphism $f \colon A \to B$ in $\mathcal{K}_0$, and $(n - 1)$-simplices $h_0, h_1, \dots, h_{n} \in \mathcal{K}(B,K)_{n - 1}$ such that $h_i \cdot f = d_i(g)$ for all ${i \in \{0, 1, \dots, n\}}$, and it returns a morphism $f' \colon K \to L$ in $\mathcal{K}_0$ and an $n$\nobreakdash-simplex ${g' \in \mathcal{K}(B,L)_n}$ such that $f' \cdot g = g' \cdot f \in \mathcal{K}(A,K)_n$ and ${f' \cdot h_i = d_i(g')}$ for all ${i \in \{0,1, \dots, n\}}$ with the following universal property: For each other such compatible pair consisting of a morphism $r \colon K \to D$ in $\mathcal{K}_0$ and an $n$-simplex ${s \in \mathcal{K}(B,D)_n}$ there exists a unique morphism $p \colon L \to D$ in $\mathcal{K}_0$ such that $p \cdot f' = r$ and $p \cdot g' = s$.
\end{enumerate}
\endexamples

\section{Factorization lemma}

In this section we prove a technical result (Lemma \ref{lem:difficultlemma}) about existence of a certain factorization in the category of morphisms that will be used in our enriched small object argument. The lemma has an easier proof under stronger assumptions, however, we do not assume these in our enriched small object argument, and thus we first present the more difficult proof of the lemma. Afterwards in the proof of Lemma \ref{lem:simplifiedlemma} we illustrate the easier proof under the stronger assumptions.
\lem\label{lem:difficultlemma}
Suppose that $\mathcal{V}$ is a cosmos, $\mathcal{K}$ is a copowered $\mathcal{V}$-category that has pushouts and transfinite composites, $g \colon X \to Y$ is a morphism in $\mathcal{K}_0$, $u \colon U \to V$ is a morphism in $\mathcal{V}_0$, $\lambda$ is a regular cardinal, $D \colon \lambda \to \mathcal{K}_0$ is a $\lambda$-sequence, $(\varphi_\alpha \colon D\alpha \to L)_{\alpha < \lambda}$ is a cocone in $\mathcal{K}_0$, the colimit of $D$ in $\mathcal{K}$ is preserved by the following three functors:
\begin{equation}\label{eqn:threefunctors}\mathcal{V}_0\big(V,\mathcal{K}(X,-)\big), \mathcal{V}_0\big(U,\mathcal{K}(X,-)\big), \mathcal{V}_0\big(U,\mathcal{K}(Y,-)\big) \colon \mathcal{K}_0 \to \mathbf{Set},\end{equation}
and $m \colon \mathrm{colim}_{\alpha < \lambda} D\alpha \to L$ is the morphism in $\mathcal{K}_0$ induced by the universal property of the colimit. Furthermore, suppose that we have the following commutative square.
\begin{equation}\label{cmd:firstassum}
\begin{tikzcd}
	{\mathrm{dom}(u \mathrel{\Box} g)} & {\operatorname{colim}_{\alpha < \lambda} D\alpha} \\
	{V \odot Y} & L
	\arrow["\ell", from=1-1, to=1-2]
	\arrow["{u \mathrel{\Box} g}"', from=1-1, to=2-1]
	\arrow["m", from=1-2, to=2-2]
	\arrow["h"', from=2-1, to=2-2]
\end{tikzcd}
\end{equation}
Then there exists an ordinal $\xi < \lambda$ and a morphism $\ell' \colon \mathrm{dom}(u \mathrel{\Box} g) \to D\xi$ such that the following diagram commutes.
\begin{equation}\label{cmd:preslemmaconclusion}
\begin{tikzcd}
	{\mathrm{dom}(u \mathrel{\Box} g)} && {\operatorname{colim}_{\alpha < \lambda} D\alpha} \\
	& {D\xi} \\
	{V \odot Y} && L
	\arrow["\ell", from=1-1, to=1-3]
	\arrow["{\ell'}"', from=1-1, to=2-2]
	\arrow["{u \mathrel{\Box} g}"', from=1-1, to=3-1]
	\arrow["m", from=1-3, to=3-3]
	\arrow["{\iota_\xi}"', from=2-2, to=1-3]
	\arrow["{\varphi_\xi}"', from=2-2, to=3-3]
	\arrow["h"', from=3-1, to=3-3]
\end{tikzcd}
\end{equation}
\endlem
\pf
Denote by $p_1$ and $p_2$ the following projections.
\begin{equation}\label{cmd:p}
\begin{tikzcd}
	{\mathcal{K}(Y,\operatorname{colim}_{\alpha < \lambda} D\alpha)} \\
	& {\mathrm{Sq}(g,m)} & {\mathcal{K}(X,\operatorname{colim}_{\alpha < \lambda} D\alpha)} \\
	& {\mathcal{K}(Y,L)} & {\mathcal{K}(X,L)}
	\arrow["{p_1}", from=2-2, to=2-3]
	\arrow["{p_2}"', from=2-2, to=3-2]
	\arrow["{\mathcal{K}(X,m)}", from=2-3, to=3-3]
	\arrow["{\mathcal{K}(g,L)}"', from=3-2, to=3-3]
	\arrow["{\mathcal{K}(Y,m)}"', from=1-1, to=3-2]
	\arrow["{\mathcal{K}(g,\operatorname{colim}_{\alpha < \lambda} D\alpha)}", from=1-1, to=2-3]
	\arrow["{\langle g,m \rangle}"{description}, dashed, from=1-1, to=2-2]
\end{tikzcd}
\end{equation}
Let
\begin{equation}\label{cmd:assum}
\begin{tikzcd}
	U & {\mathcal{K}(Y,\operatorname{colim}_{\alpha < \lambda} D\alpha)} \\
	V & {\mathrm{Sq}(g,m)}
	\arrow["{(\ell \cdot i_1)_*}", from=1-1, to=1-2]
	\arrow["u"', from=1-1, to=2-1]
	\arrow["{\langle g,m \rangle}", from=1-2, to=2-2]
	\arrow["{((\ell \cdot i_2)_*,h_*)}"', from=2-1, to=2-2]
\end{tikzcd}
\end{equation}
be the adjoint transpose (recall Remark \ref{rem:enabadj}) of the commutative square (\ref{cmd:firstassum}). To simplify notation, define $w := (\ell \cdot i_1)_*$ and $v := ((\ell \cdot i_2)_*,h_*)$. By assumption, we know that
$$\Big(\mathcal{V}_0\big(V,\mathcal{K}(X,\iota_\beta)\big) \colon \mathcal{V}_0\big(V,\mathcal{K}(X,D\beta)\big) \to \mathcal{V}_0\big(V,\mathcal{K}(X,\operatorname{colim}_{\alpha < \lambda} D\alpha)\big)\Big)_{\beta < \lambda}$$
is a colimit in $\mathbf{Set}$, and thus, by the characterization of directed colimits in $\mathbf{Set}$, there exists $\xi < \lambda$ and a morphism $v'' \colon V \to \mathcal{K}(X, D\xi)$ such that the following square commutes.
\begin{equation}\label{eqn:v}
\begin{tikzcd}
	V & {\mathcal{K}(X,D\xi)} \\
	{\mathrm{Sq}(g,m)} & {\mathcal{K}(X,\operatorname{colim}_{\alpha < \lambda} D\alpha)}
	\arrow["{v''}", from=1-1, to=1-2]
	\arrow["v"', from=1-1, to=2-1]
	\arrow["{\mathcal{K}(X,\iota_\xi)}", from=1-2, to=2-2]
	\arrow["{p_1}"', from=2-1, to=2-2]
\end{tikzcd}
\end{equation}
Furthermore, again by assumption, we know that
$$\Big(\mathcal{V}_0\big(U,\mathcal{K}(Y,\iota_\beta)\big) \colon \mathcal{V}_0\big(U,\mathcal{K}(Y,D\beta)\big) \to \mathcal{V}_0\big(U,\mathcal{K}(Y,\operatorname{colim}_{\alpha < \lambda} D\alpha)\big)\Big)_{\beta < \lambda}$$
is a colimit in $\mathbf{Set}$, and thus, by the characterization of directed colimits in $\mathbf{Set}$, there exists an ordinal $\xi' < \lambda$ and a morphism $w' \colon U \to \mathcal{K}(Y, D\xi')$ such that $w = \mathcal{K}(Y, \iota_{\xi'}) \cdot w'$. By using $\min\{\xi, \xi'\} \to \max\{\xi, \xi'\}$ we can assume that $\xi' = \xi$, and thus the following triangle commutes.
\begin{equation}\label{eqn:w}
\begin{tikzcd}
	& {\mathcal{K}(Y,D\xi)} \\
	U & {\mathcal{K}(Y,\operatorname{colim}_{\alpha < \lambda} D\alpha)}
	\arrow["{\mathcal{K}(Y,\iota_{\xi})}", from=1-2, to=2-2]
	\arrow["{w'}", from=2-1, to=1-2]
	\arrow["w"', from=2-1, to=2-2]
\end{tikzcd}
\end{equation}
Using the fact that $\mathrm{Sq}(g,\varphi_\xi)$ is a pullback we get a unique morphism $v'$ making the two triangles in the following diagram commute.
\begin{equation}\label{cmd:v'}
\begin{tikzcd}
	{V} \\
	& {\mathrm{Sq}(g,\varphi_\xi)} & {\mathcal{K}(X,D\xi)} \\
	& {\mathcal{K}(Y,L)} & {\mathcal{K}(X,L)}
	\arrow["{\pi_1}"', from=2-2, to=2-3]
	\arrow["{\pi_2}"', from=2-2, to=3-2]
	\arrow["{\mathcal{K}(X,\varphi_\xi)}", from=2-3, to=3-3]
	\arrow["{\mathcal{K}(g,L)}"', from=3-2, to=3-3]
	\arrow["{p_2 \cdot v}"', curve={height=6pt}, from=1-1, to=3-2]
	\arrow["{v''}", from=1-1, to=2-3]
	\arrow["{v'}"{description}, dashed, from=1-1, to=2-2]
\end{tikzcd}
\end{equation}
Indeed, $$\mathcal{K}(g,L) \cdot p_2 \cdot v \stackrel{(\ref{cmd:p})}{=} \mathcal{K}(X,m) \cdot p_1 \cdot v \stackrel{(\ref{eqn:v})}{=} \mathcal{K}(X,m) \cdot \mathcal{K}(X, \iota_\xi) \cdot v'' = \mathcal{K}(X,\varphi_\xi) \cdot v'',$$
where the last equality follows from the definition of $m$. We will now show that the following square commutes.
\begin{equation}\label{cmd:difficultsquare}
\begin{tikzcd}
	U & {\mathcal{K}(Y,D\xi)} \\
	V & {\mathrm{Sq}(g,\varphi_\xi)}
	\arrow["{w'}", from=1-1, to=1-2]
	\arrow["u"', from=1-1, to=2-1]
	\arrow["{\langle g, \varphi_\xi \rangle}", from=1-2, to=2-2]
	\arrow["{v'}"', from=2-1, to=2-2]
\end{tikzcd}
\end{equation}
Indeed, note that
\begin{align*}
\mathcal{K}(X,\iota_\xi) \cdot \pi_1 \cdot v' \cdot u &\stackrel{\mathclap{(\ref{cmd:v'})}}{=} \mathcal{K}(X,\iota_\xi) \cdot v'' \cdot u \\
&\stackrel{\mathclap{(\ref{eqn:v})}}{=} p_1 \cdot v \cdot u \\ 
&\stackrel{\mathclap{(\ref{cmd:assum})}}{=} p_1 \cdot \langle g,m \rangle \cdot w \\ 
&\stackrel{\mathclap{(\ref{cmd:p})}}{=} \mathcal{K}(g,\operatorname{colim}_{\alpha < \lambda} D\alpha) \cdot w \\ 
&\stackrel{\mathclap{(\ref{eqn:w})}}{=} \mathcal{K}(g,\operatorname{colim}_{\alpha < \lambda} D\alpha) \cdot \mathcal{K}(Y,\iota_\xi) \cdot w' \\ 
&= \mathcal{K}(X,\iota_\xi) \cdot \mathcal{K}(g, D\xi) \cdot w' \\
&\stackrel{\mathclap{(\ref{cmd:inducede})}}{=} \mathcal{K}(X,\iota_\xi) \cdot \pi_1 \cdot \langle g, \varphi_\xi \rangle \cdot w'.
\end{align*}
Therefore, since by assumption
$$\Big(\mathcal{V}_0\big(U,\mathcal{K}(X,\iota_\beta)\big) \colon \mathcal{V}_0\big(U,\mathcal{K}(X,D\beta)\big) \to \mathcal{V}_0\big(U,\mathcal{K}(X,\operatorname{colim}_{\alpha < \lambda} D\alpha)\big)\Big)_{\beta < \lambda}$$
is a colimit in $\mathbf{Set}$, using the characterization of directed colimits in $\mathbf{Set}$, we get that there exists $\tilde{\xi}$ such that $\lambda > \tilde{\xi} \geq \xi$, 
\begin{equation}\label{eqn:xitotildexi}
\mathcal{K}(X,D(\xi \to \tilde{\xi})) \cdot \pi_1 \cdot v' \cdot u = \mathcal{K}(X,D(\xi \to \tilde{\xi})) \cdot \pi_1 \cdot \langle g, \varphi_\xi \rangle \cdot w',
\end{equation}
and morphisms $v''_{\tilde{\xi}} = \mathcal{K}(X,D(\xi \to \tilde{\xi})) \cdot v''$, $w'_{\tilde{\xi}} = \mathcal{K}(Y,D(\xi \to \tilde{\xi})) \cdot w',$ and $v'_{\tilde{\xi}}$ satisfying the analogues of (\ref{eqn:v}), (\ref{eqn:w}), and (\ref{cmd:v'}) for $\tilde{\xi}$, respectively. Thus we obtain the following chain of equalities:
\begin{align*}
\pi_{1,\tilde{\xi}} \cdot v'_{\tilde{\xi}} \cdot u &\stackrel{\mathclap{(\tilde{\ref{cmd:v'}})}}{=}  v''_{\tilde{\xi}} \cdot u\\
&= \mathcal{K}\big(X,D(\xi \to \tilde{\xi})\big) \cdot v'' \cdot u \\
&\stackrel{\mathclap{(\ref{cmd:v'})}}{=} \mathcal{K}\big(X,D(\xi \to \tilde{\xi})\big) \cdot \pi_1 \cdot v' \cdot u \\
&\stackrel{\mathclap{(\ref{eqn:xitotildexi})}}{=} \mathcal{K}\big(X,D(\xi \to \tilde{\xi})\big) \cdot \pi_1 \cdot \langle g, \varphi_\xi \rangle \cdot w' \\ 
&\stackrel{\mathclap{(\ref{cmd:inducede})}}{=} \mathcal{K}\big(X,D(\xi \to \tilde{\xi})\big) \cdot \mathcal{K}(g,D\xi) \cdot w' \\ 
&= \mathcal{K}(g,D\tilde{\xi}) \cdot \mathcal{K}\big(Y,D(\xi \to \tilde{\xi})\big) \cdot w' \\ 
&= \mathcal{K}(g,D\tilde{\xi}) \cdot w'_{\tilde{\xi}} \\
&\stackrel{\mathclap{(\ref{cmd:inducede})}}{=} \pi_{1,\tilde{\xi}} \cdot \langle g, \varphi_{\tilde{\xi}} \rangle \cdot w'_{\tilde{\xi}}.
\end{align*}
Now we can replace the previous $\xi$ by $\tilde{\xi}$ (while still denoting it $\xi$) and we obtain $$\pi_1 \cdot v' \cdot u = \pi_1 \cdot \langle g, \varphi_\xi \rangle \cdot w'.$$
Furthermore,
\begin{align*}
\pi_2 \cdot v' \cdot u &\stackrel{\mathclap{(\ref{cmd:v'})}}{=} p_2 \cdot v \cdot u \\
&\stackrel{\mathclap{(\ref{cmd:assum})}}{=} p_2 \cdot \langle g,m \rangle \cdot w \\
&\stackrel{\mathclap{(\ref{cmd:p})}}{=} \mathcal{K}(Y,m) \cdot w \\
&\stackrel{\mathclap{(\ref{eqn:w})}}{=} \mathcal{K}(Y,m) \cdot \mathcal{K}(Y,\iota_\xi) \cdot w' \\
&= \mathcal{K}(Y,\varphi_\xi) \cdot w' \\
&\stackrel{\mathclap{(\ref{cmd:inducede})}}{=} \pi_2 \cdot \langle g, \varphi_\xi \rangle \cdot w'.
\end{align*}
Since $(\pi_1, \pi_2)$ is a pullback, we obtain that the square (\ref{cmd:difficultsquare}) commutes. Consider the adjoint transpose of (\ref{cmd:difficultsquare}), which gives us the following commutative square.
\begin{equation}\label{cmd:bottomleftsq}
\begin{tikzcd}[column sep = 4em]
	{\mathrm{dom}(u \mathrel{\Box} g)} & {D\xi} \\
	{V \odot Y} & L
	\arrow["{((w')^*,(\pi_1 \cdot v')^*)}", from=1-1, to=1-2]
	\arrow["{u \mathrel{\Box} g}"', from=1-1, to=2-1]
	\arrow["{\varphi_\xi}", from=1-2, to=2-2]
	\arrow["{(\pi_2 \cdot v')^*}"', from=2-1, to=2-2]
\end{tikzcd}
\end{equation}
Note that $(\pi_2 \cdot v')^* = (p_2 \cdot v)^* = (h_*)^* = h$, and thus we can define $\ell' := ((w')^*,(\pi_1 \cdot v')^*)$ and then the commutative square (\ref{cmd:bottomleftsq}) demonstrates that the bottom-left square in the diagram (\ref{cmd:preslemmaconclusion}) commutes. Furthermore,
$$\iota_\xi \cdot \ell' = \iota_\xi \cdot ((w')^*,(\pi_1 \cdot v')^*) \stackrel{(\ref{cmd:v'})}{=} \iota_\xi \cdot ((w')^*,(v'')^*) \stackrel{(\ref{eqn:v}) \& (\ref{eqn:w})}{=} (w^*, (p_1 \cdot v)^*) = (\ell \cdot i_1, \ell \cdot i_2) = \ell,$$
and hence the morphism $\ell'$ makes the diagram (\ref{cmd:preslemmaconclusion}) commute.
\epf

A simpler proof (which was suggested to the author by the anonymous referee) of Lemma \ref{lem:difficultlemma} can be given under the following additional assumption: On top of assuming preservation by the three functors (\ref{eqn:threefunctors}), assume also preservation by the functor $\mathcal{V}_0\big(V,\mathcal{K}(Y,-)\big) \colon \mathcal{K}_0 \to \mathbf{Set}$. In fact, with this additional assumption one can obtain an even stronger conclusion (which, however, is not needed for our small object argument). For simplicity of phrasing, let us illustrate the simpler approach by proving the following lemma that uses the language of presentability.

\lem\label{lem:simplifiedlemma}
Suppose that $\mathcal{V}$ is a cosmos, $\kappa$ is a regular cardinal, $\mathcal{K}$ is a copowered $\mathcal{V}$\nobreakdash-category that has pushouts and $\kappa$-directed colimits, $g \colon X \to Y$ is a morphism in $\mathcal{K}_0$, $u \colon U \to V$ is a morphism in $\mathcal{V}_0$, the objects $X$ and $Y$ are $\kappa$-presentable in the enriched sense \mbox{\cite[2.1]{Kel82}}, and the objects $U$ and $V$ are $\kappa$-presentable in the unenriched sense. Then $u \mathrel{\Box} g$ is $\kappa$\nobreakdash-presentable in $[\mathbf{2},\mathcal{K}]_0$ in the unenriched sense.
\endlem
\pf
Suppose that $D \colon I \to [\mathbf{2},\mathcal{K}]_0$ is a $\kappa$-directed diagram. For each $i \in I$, we will denote $D(i)$ by $f_i \colon A_i \to B_i$. We wish to show that
\[\operatorname{colim}_{i \in I} [\mathbf{2},\mathcal{K}]_0(u \mathrel{\Box} g, f_i) \cong [\mathbf{2},\mathcal{K}]_0(u \mathrel{\Box} g, \operatorname{colim}_{i \in I} f_i)\]
in $\mathbf{Set}$. By Remark \ref{rem:enabadj}, this is equivalent to showing that
\[\operatorname{colim}_{i \in I} [\mathbf{2},\mathcal{V}]_0\big(u, \langle g, f_i\rangle\big) \cong [\mathbf{2},\mathcal{V}]_0\big(u, \langle g, \operatorname{colim}_{i \in I} f_i\rangle\big)\]
in $\mathbf{Set}$. By definition, a morphism $u \to \langle g, \operatorname{colim}_{i \in I} f_i\rangle$ consists of two components that make the obvious square commute, and thus we obtain the following pullback in $\mathbf{Set}$.
\[\begin{tikzcd}
	{[\mathbf{2},\mathcal{V}]_0\big(u,\langle g, \operatorname{colim}_{i \in I} f_i\rangle\big)} & {\mathcal{V}_0\big(U,\mathcal{K}(Y, \operatorname{colim}_{i \in I} A_i)\big)} \\
	{\mathcal{V}_0\big(V,\operatorname{Sq}(g, \operatorname{colim}_{i \in I} f_i)\big)} & {\mathcal{V}_0\big(U,\operatorname{Sq}(g, \operatorname{colim}_{i \in I} f_i)\big)}
	\arrow["{\pi_1}", from=1-1, to=1-2]
	\arrow["{\pi_2}"', from=1-1, to=2-1]
	\arrow["{\mathcal{V}_0(U,\langle g, \operatorname{colim}_{i \in I} f_i\rangle)}", from=1-2, to=2-2]
	\arrow["{\mathcal{V}_0(u,\operatorname{Sq}(g, \operatorname{colim}_{i \in I} f_i))}"', shift right=2, draw=none, from=2-1, to=2-2]
	\arrow[from=2-1, to=2-2]
\end{tikzcd}\]
Using the fact that $\kappa$-directed colimits commute with pullbacks in $\mathbf{Set}$ we see that it suffices to show the corresponding isomorphisms in $\mathbf{Set}$ for the three other objects in the pullback square. Let us begin with the object $\mathcal{V}_0\big(U,\mathcal{K}(Y, \operatorname{colim}_{i \in I} A_i)\big)$, since that case is the easiest. By assumption, $Y$ is $\kappa$-presentable in the enriched sense and $U$ is $\kappa$-presentable in the unenriched sense. Therefore,
\[\mathcal{V}_0\big(U,\mathcal{K}(Y, \operatorname{colim}_{i \in I} A_i)\big) \cong \mathcal{V}_0\big(U,\operatorname{colim}_{i \in I} \mathcal{K}(Y, A_i)\big) \cong \operatorname{colim}_{i \in I} \mathcal{V}_0\big(U,\mathcal{K}(Y, A_i)\big).\]
The two remaining cases are analogous to each other, and thus we will describe only the case of the object $\mathcal{V}_0\big(V,\operatorname{Sq}(g, \operatorname{colim}_{i \in I} f_i)\big)$. By definition of $\operatorname{Sq}(g, \operatorname{colim}_{i \in I} f_i)$, we have the following pullback
\[\begin{tikzcd}
	{\operatorname{Sq}(g, \operatorname{colim}_{i \in I} f_i)} & {\mathcal{K}(X, \operatorname{colim}_{i \in I} A_i)} \\
	{\mathcal{K}(Y, \operatorname{colim}_{i \in I} B_i)} & {\mathcal{K}(X, \operatorname{colim}_{i \in I} B_i)}
	\arrow["{p_1}", from=1-1, to=1-2]
	\arrow["{p_2}"', from=1-1, to=2-1]
	\arrow["{\mathcal{K}(X, \operatorname{colim}_{i \in I} f_i)}", from=1-2, to=2-2]
	\arrow["{\mathcal{K}(g, \operatorname{colim}_{i \in I} B_i)}"', shift right, draw=none, from=2-1, to=2-2]
	\arrow[from=2-1, to=2-2]
\end{tikzcd}\]
and thus, since $\mathcal{V}_0(V,-) \colon \mathcal{V}_0 \to \mathbf{Set}$ preserves limits, we obtain the following pullback in $\mathbf{Set}$.
\[\begin{tikzcd}
	{\mathcal{V}_0\big(V,\operatorname{Sq}(g, \operatorname{colim}_{i \in I} f_i)\big)} & {\mathcal{V}_0\big(V,\mathcal{K}(X, \operatorname{colim}_{i \in I} A_i)\big)} \\
	{\mathcal{V}_0\big(V,\mathcal{K}(Y, \operatorname{colim}_{i \in I} B_i)\big)} & {\mathcal{V}_0\big(V,\mathcal{K}(X, \operatorname{colim}_{i \in I} B_i)\big)}
	\arrow["{\mathcal{V}_0(V,p_1)}", from=1-1, to=1-2]
	\arrow["{\mathcal{V}_0(V,p_2)}"', from=1-1, to=2-1]
	\arrow["{\mathcal{V}_0(V,\mathcal{K}(X, \operatorname{colim}_{i \in I} f_i))}", from=1-2, to=2-2]
	\arrow["{\mathcal{V}_0(V,\mathcal{K}(g, \operatorname{colim}_{i \in I} B_i))}"', shift right=2, draw=none, from=2-1, to=2-2]
	\arrow[from=2-1, to=2-2]
\end{tikzcd}\]
Therefore, in order to conclude that $\operatorname{colim}_{i \in I} \mathcal{V}_0\big(V,\operatorname{Sq}(g,f_i)\big) \cong \mathcal{V}_0\big(V,\operatorname{Sq}(g, \operatorname{colim}_{i \in I} f_i)\big)$, it suffices to show the corresponding isomorphisms for the three other objects in the pullback, since $\kappa$-directed colimits commute with pullbacks in $\mathbf{Set}$. Indeed, by the $\kappa$-presentability assumptions on $X$, $Y$, $V$ we obtain the following isomorphisms
\begin{align*}
\mathcal{V}_0\big(V,\mathcal{K}(X, \operatorname{colim}_{i \in I} A_i)\big) &\cong \mathcal{V}_0\big(V,\operatorname{colim}_{i \in I} \mathcal{K}(X, A_i)\big) \cong \operatorname{colim}_{i \in I} \mathcal{V}_0\big(V,\mathcal{K}(X, A_i)\big), \\
\mathcal{V}_0\big(V,\mathcal{K}(X, \operatorname{colim}_{i \in I} B_i)\big) &\cong \mathcal{V}_0\big(V,\operatorname{colim}_{i \in I} \mathcal{K}(X, B_i)\big) \cong \operatorname{colim}_{i \in I} \mathcal{V}_0\big(V,\mathcal{K}(X, B_i)\big), \\
\mathcal{V}_0\big(V,\mathcal{K}(Y, \operatorname{colim}_{i \in I} B_i)\big) &\cong \mathcal{V}_0\big(V,\operatorname{colim}_{i \in I} \mathcal{K}(Y, B_i)\big) \cong \operatorname{colim}_{i \in I} \mathcal{V}_0\big(V,\mathcal{K}(Y, B_i)\big).
\end{align*}
Note that the last line corresponds to the aforementioned preservation assumption that's not assumed in Lemma \ref{lem:difficultlemma}.
\epf

\section{Enriched small object argument}

Here we finally state and prove the small object argument. We prove the small object argument in the special case when $\mathscr{F}$ is a cofibrantly generated weak factorization system. Thus, in this final section we assume that $\mathcal{V}$ is a cosmos, $\mathcal{K}$ is a copowered $\mathcal{V}$-category that has pushouts and transfinite composites, $\mathcal{I}$ is a set of morphisms in $\mathcal{K}_0$, $\mathcal{J}$ is a set of morphisms in $\mathcal{V}_0$, and $\mathscr{F} = (\mathscr{L}, \mathscr{R})$ is a weak factorization system on $\mathcal{V}_0$ that's cofibrantly generated by $\mathcal{J}$. Before stating and proving the theorem we give two definitions concerning smallness and relative cell complexes.
\definition\label{def:newsmallness}
Suppose that $\mathcal{S}$ is a class of morphisms in $\mathcal{K}_0$, $U$ and $V$ are objects in $\mathcal{V}$, $X$ and $Y$ are objects in $\mathcal{K}$, and $\kappa$ is a cardinal. We say that $(U,V,X,Y)$ is \emph{$\kappa$-small relative to $\mathcal{S}$} if for all $\kappa$-filtered ordinals $\lambda$ and all $\lambda$-sequences $D \colon \lambda \to \mathcal{K}_0$ such that for each $\beta + 1 < \lambda$ the morphism $D(\beta \to \beta + 1) \colon D\beta \to D(\beta + 1)$ belongs to $\mathcal{S}$, it is the case that the colimit of $D$ in $\mathcal{K}$ is preserved by the following three functors:
$$\mathcal{V}_0\big(V,\mathcal{K}(X,-)\big), \mathcal{V}_0\big(U,\mathcal{K}(X,-)\big), \mathcal{V}_0\big(U,\mathcal{K}(Y,-)\big) \colon \mathcal{K}_0 \to \mathbf{Set}.$$
Furthermore, we say that $(U,V,X,Y)$ is \emph{small relative to $\mathcal{S}$} if there exists a cardinal $\kappa$ such that $(U,V,X,Y)$ is $\kappa$-small relative to $\mathcal{S}$.
\enddefinition
\example
If $\mathcal{V} = \mathbf{Set}$, then smallness of $(\emptyset,1,X,Y)$ relative to $\mathcal{S}$ from Definition \ref{def:newsmallness} is the same as the ordinary smallness of $X$ relative to $\mathcal{S}$ \cite[Definition 2.1.3]{Hov99}. Subsequently, this means that the enriched small object argument (Theorem \ref{thm:soa}) that we prove below includes the ordinary small object argument as a special case. We will explain this later on in more detail in Example \ref{ex:enrichedsoafinal}.(\ref{ex:insSet}).
\endexample
\rem
Intuitively, the smallness condition on $(U,V,X,Y)$ in Definition \ref{def:newsmallness} involves an enriched smallness condition on $X$ and $Y$ and an unenriched smallness condition on $U$ and $V$. The presentable case (Remark \ref{rem:replacesmallnessbypres}) will be particularly illustrative for this intuition.
\endrem
\definition
By $(\mathcal{I} \mathrel{\mathrm{rel}} \mathcal{J})\text{-}\mathrm{cell}$ we denote the class of all transfinite composites of copowered pushouts of morphisms from $\mathcal{I}$ relative to $\mathcal{J}$.
\enddefinition
\thm\label{thm:soa}
Suppose that $\mathcal{V}$ is a cosmos, $\mathcal{K}$ is a copowered $\mathcal{V}$-category that has pushouts and transfinite composites, $\mathcal{I}$ is a set of morphisms in $\mathcal{K}_0$, $\mathcal{J}$ is a set of morphisms in $\mathcal{V}_0$, $\mathscr{F} = (\mathscr{L}, \mathscr{R})$ is a weak factorization system on $\mathcal{V}_0$ that's cofibrantly generated by $\mathcal{J}$, $\mathscr{L}$ is stable under corners, and for each morphism $u \colon U \to V$ in $\mathcal{J}$ and each morphism $g \colon X \to Y$ in $\mathcal{I}$ it is the case that $(U,V,X,Y)$ is small relative to $(\mathcal{I} \mathrel{\mathrm{rel}} \mathcal{J})\text{-}\mathrm{cell}$.\par
Then for each morphism $f \colon K \to L$ in $\mathcal{K}_0$ there exists a factorization $f = m \cdot e$ such that $e$ and $m$ are morphisms in $\mathcal{K}_0$ satisfying $e \in (\mathcal{I} \mathrel{\mathrm{rel}} \mathcal{J})\text{-}\mathrm{cell}$ and $m \in \mathcal{I}^{\overset{\mathscr{F}}{\pitchfork}}$.
\endthm
\pf
Let $\kappa$ be a cardinal such that for each morphism $u \colon U \to V$ in $\mathcal{J}$ and each morphism $g \colon X \to Y$ in $\mathcal{I}$ it is the case that $(U,V,X,Y)$ is $\kappa$-small relative to the class $(\mathcal{I} \mathrel{\mathrm{rel}} \mathcal{J})\text{-}\mathrm{cell}$. Furthermore, let $\lambda$ be a $\kappa$-filtered regular cardinal such that the cardinality of $\mathcal{J}$ is less than $\lambda$.\par
We will transfinitely define a diagram $D \colon \lambda \to \mathcal{K}_0$ such that $D0 = K$, the morphism ${D(\alpha \to \alpha + 1)}$ is a transfinite composition of copowered pushouts of morphisms from $\mathcal{I}$ relative to $\mathcal{J}$ for each ordinal $\alpha < \lambda$, and $D\alpha := \mathrm{colim}_{\zeta < \alpha} D\zeta$ for each limit ordinal $\alpha < \lambda$. Then we will define $e$ to be the injection ${\iota_0 \colon D0 \to \mathrm{colim}_{\alpha < \lambda} D\alpha}$. Furthermore, we are going to construct a cocone $(\varphi_\alpha \colon D\alpha \to L)_{\alpha < \lambda}$ in $\mathcal{K}_0$, and then we can define the morphism ${m \colon \mathrm{colim}_{\alpha < \lambda} D\alpha \to L}$ to be the induced morphism obtained from the universal property of the colimit.\newline
Base Case: Define $D0 := K$ and $\varphi_0 := f$.\newline
Limit Step: Suppose that $\alpha < \lambda$ is a limit ordinal such that we've already performed the construction for all $\zeta < \alpha$. Then define $D\alpha := \operatorname{colim}_{\zeta < \alpha} D\zeta$ and let $\varphi_{\alpha} \colon D\alpha \to L$ be the morphism induced by the cocone $(\varphi_{\zeta} \colon D\zeta \to L)_{\zeta < \alpha}$.\newline
Successor Step: Suppose that $\alpha < \lambda$ is an ordinal such that we've already performed the construction for each $\zeta \leq \alpha$. Enumerate the set $\mathcal{J}$ by its cardinality, which means that $\mathcal{J} = \{u_{\alpha'} \colon U_{\alpha'} \to V_{\alpha'} \mid \alpha' < \alpha_{\mathcal{J}}\}$, where $\alpha_{\mathcal{J}}$ is the cardinality of $\mathcal{J}$ viewed as an ordinal number. Let $F \colon \mathbf{Ord} \to \mathbf{Ord}$ be the ordinal function defined via transfinite induction as follows: $F(0) := 0$, if $F$ is defined on $\gamma$, then
\begin{equation*}
F(\gamma + 1) :=
\begin{cases}
          F(\gamma) + 1 \quad &\text{if} \, F(\gamma) + 1 < \alpha_{\mathcal{J}},\\
          0 \quad &\text{if} \, F(\gamma) + 1 \geq \alpha_{\mathcal{J}},\\
\end{cases}
\end{equation*}
and finally if $\gamma$ is a limit ordinal, then we consider two cases:
\begin{itemize}
\item if $\alpha_\mathcal{J}$ is a finite ordinal, then $F(\gamma) := 0$, and
\item if $\alpha_{\mathcal{J}}$ is a limit ordinal, then
\begin{equation*}
F(\gamma) :=
\begin{cases}
          \sup_{\delta < \gamma}F(\delta) \quad &\text{if} \, \sup_{\delta < \gamma}F(\delta) < \alpha_\mathcal{J},\\
          0 \quad &\text{if} \, \sup_{\delta < \gamma}F(\delta) \geq \alpha_\mathcal{J}.\\
\end{cases}
\end{equation*}
\end{itemize}
Briefly, we can think of $F$ as being defined via the modulo operator in the following way $F := -\! \mod \alpha_{\mathcal{J}}$. Consider the set $J_\alpha$ of all triples $(g,v,w)$ of morphisms such that $g \in \mathcal{I}$ and that make the following square commute.
\[\begin{tikzcd}
	{\mathrm{dom}(u_{F(\alpha)} \mathrel{\Box} g)} & {D\alpha} \\
	{V_{F(\alpha)} \odot \operatorname{cod} g} & L
	\arrow["w", from=1-1, to=1-2]
	\arrow["{u_{F(\alpha)} \mathrel{\Box} g}"', from=1-1, to=2-1]
	\arrow["{\varphi_\alpha}", from=1-2, to=2-2]
	\arrow["v"', from=2-1, to=2-2]
\end{tikzcd}\]
Now we are going to use transfinite induction again. Well-order the set $J_\alpha$, which means that ${J_\alpha = \{(g_\gamma, v_\gamma, w_\gamma) \mid \gamma < \gamma_\alpha\}}$, where $\gamma_\alpha$ is an ordinal. We define $D_0(\alpha + 1) := D\alpha$ and ${\varphi_{\alpha + 1, 0} := \varphi_\alpha}$. If $\beta \leq \gamma_\alpha$ is a limit ordinal such that we've already performed the construction for all the ordinals less than $\beta$, then define $D_\beta(\alpha + 1) := \operatorname{colim}_{\delta < \beta} D_\delta(\alpha + 1)$ and define $\varphi_{\alpha + 1,\beta} \colon D_{\beta}(\alpha + 1) \to L$ to be the morphism that's induced by the cocone $(\varphi_{\alpha + 1,\delta} \colon D_{\delta}(\alpha + 1) \to L)_{\delta < \beta}$. Finally, suppose that $\beta < \gamma_\alpha$ is an ordinal such that we've already performed the construction for all the ordinals less or equal to $\beta$. Denote by $f'_{0,\beta} \colon D\alpha \to D_{\beta}(\alpha + 1)$ the morphism coming from the construction in the previous steps. Define $D_{\beta + 1}(\alpha + 1)$ to be following pushout in $\mathcal{K}$.
\[\begin{tikzcd}
	{\mathrm{dom}(u_{F(\alpha)} \mathrel{\Box} g_\beta)} & {D_\beta(\alpha + 1)} \\
	{V_{F(\alpha)} \odot \operatorname{cod} g_\beta} & {D_{\beta + 1}(\alpha + 1)}
	\arrow["{f'_{0,\beta} \cdot w_\beta}", from=1-1, to=1-2]
	\arrow["{u_{F(\alpha)} \mathrel{\Box} g_\beta}"', from=1-1, to=2-1]
	\arrow["{f'_{\beta,\beta + 1}}", from=1-2, to=2-2]
	\arrow["{h'_\beta}"', from=2-1, to=2-2]
\end{tikzcd}\]
Furthermore, define $\varphi_{\alpha + 1,\beta + 1} \colon D_{\beta + 1}(\alpha + 1) \to L$ to be the morphism induced by the following commutative square.
\[\begin{tikzcd}
	{\mathrm{dom}(u_{F(\alpha)} \mathrel{\Box} g_\beta)} & {D_\beta(\alpha + 1)} \\
	{V_{F(\alpha)} \odot \operatorname{cod} g_\beta} & L
	\arrow["{f'_{0,\beta} \cdot w_\beta}", from=1-1, to=1-2]
	\arrow["{u_{F(\alpha)} \mathrel{\Box} g_\beta}"', from=1-1, to=2-1]
	\arrow["{\varphi_{\alpha + 1,\beta}}", from=1-2, to=2-2]
	\arrow["{v_\beta}"', from=2-1, to=2-2]
\end{tikzcd}\]
The square indeed commutes, since $\varphi_{\alpha + 1,\beta} \cdot f'_{0,\beta} \cdot w_\beta = \varphi_{\alpha} \cdot w_\beta = v_\beta \cdot (u_{F(\alpha)} \mathrel{\Box} g_\beta)$. The ``inner'' transfinite construction is now finished. Define ${D(\alpha \to \alpha + 1)}$ to be the transfinite composition of all the morphisms $f'_{\delta, \delta + 1}$, where $\delta < \gamma_\alpha$, and define the morphism $\varphi_{\alpha + 1} \colon D(\alpha + 1) \to L$ to be the map induced by the cocone ${(\varphi_{\alpha + 1,\beta} \colon D_{\beta}(\alpha + 1) \to L)_{\beta < \gamma_\alpha}}$.\par
Now that we've finished the construction let us show that $m \in \mathcal{I}^{\overset{\mathscr{F}}{\pitchfork}}$. Suppose that $g \colon X \to Y$ is in $\mathcal{I}$. We want to verify that $\langle g,m \rangle \in \mathcal{J}^\pitchfork$, where $\langle g,m \rangle$ is the induced morphism below.
\begin{equation*}
\begin{tikzcd}
	{\mathcal{K}(Y,\operatorname{colim}_{\alpha < \lambda} D\alpha)} \\
	& {\mathrm{Sq}(g,m)} & {\mathcal{K}(X,\operatorname{colim}_{\alpha < \lambda} D\alpha)} \\
	& {\mathcal{K}(Y,L)} & {\mathcal{K}(X,L)}
	\arrow["{p_1}", from=2-2, to=2-3]
	\arrow["{p_2}"', from=2-2, to=3-2]
	\arrow["{\mathcal{K}(X,m)}", from=2-3, to=3-3]
	\arrow["{\mathcal{K}(g,L)}"', from=3-2, to=3-3]
	\arrow["{\mathcal{K}(Y,m)}"', from=1-1, to=3-2]
	\arrow["{\mathcal{K}(g,\operatorname{colim}_{\alpha < \lambda} D\alpha)}", from=1-1, to=2-3]
	\arrow["{\langle g,m \rangle}"{description}, dashed, from=1-1, to=2-2]
\end{tikzcd}
\end{equation*}
Suppose that we have a commutative square
\begin{equation*}
\begin{tikzcd}
	{U_{\alpha'}} & {\mathcal{K}(Y,\operatorname{colim}_{\alpha < \lambda} D\alpha)} \\
	{V_{\alpha'}} & {\mathrm{Sq}(g,m)}
	\arrow["{\langle g,m \rangle}", from=1-2, to=2-2]
	\arrow["v"', from=2-1, to=2-2]
	\arrow["w", from=1-1, to=1-2]
	\arrow["{u_{\alpha'}}"', from=1-1, to=2-1]
\end{tikzcd}
\end{equation*}
in which $u_{\alpha'} \in \mathcal{J}$. We want to show that there exists a diagonal $V_{\alpha'} \to \mathcal{K}(Y,\operatorname{colim}_{\alpha < \lambda} D\alpha)$ that makes both triangles commute. By Remark \ref{rem:enabadj}, this is equivalent to showing that there exists a diagonal $V_{\alpha'} \odot Y \to \operatorname{colim}_{\alpha < \lambda} D\alpha$ that makes the two triangles in the following square commute.
\[\begin{tikzcd}[column sep = 3em]
	{\mathrm{dom}(u_{\alpha'} \mathrel{\Box} g)} & {\operatorname{colim}_{\alpha < \lambda} D\alpha} \\
	{V_{\alpha'} \odot Y} & L
	\arrow["{(w^*,(p_1 \cdot v)^*)}", from=1-1, to=1-2]
	\arrow["{u_{\alpha'} \mathrel{\Box} g}"', from=1-1, to=2-1]
	\arrow["m", from=1-2, to=2-2]
	\arrow["{(p_2 \cdot v)^*}"', from=2-1, to=2-2]
\end{tikzcd}\]
Thus, let us show the existence of such a diagonal $V_{\alpha'} \odot Y \to \operatorname{colim}_{\alpha < \lambda} D\alpha$. By Lemma \ref{lem:difficultlemma} we get that there exists $\xi < \alpha$ and $w' \colon \mathrm{dom}(u_{\alpha'} \mathrel{\Box} g) \to D\xi$ such that the following diagram commutes.
\[\begin{tikzcd}
	{\mathrm{dom}(u_{\alpha'} \mathrel{\Box} g)} && {\operatorname{colim}_{\alpha < \lambda} D\alpha} \\
	& {D\xi} \\
	{V_{\alpha'} \odot Y} && L
	\arrow["{(w^*,(p_1 \cdot v)^*)}", from=1-1, to=1-3]
	\arrow["{w'}"', from=1-1, to=2-2]
	\arrow["{u_{\alpha'} \mathrel{\Box} g}"', from=1-1, to=3-1]
	\arrow["m", from=1-3, to=3-3]
	\arrow["{\iota_\xi}"', from=2-2, to=1-3]
	\arrow["{\varphi_\xi}"', from=2-2, to=3-3]
	\arrow["{(p_2 \cdot v)^*}"', from=3-1, to=3-3]
\end{tikzcd}\]
Furthermore, since the cardinality of $\mathcal{J}$ is less than $\lambda$, there exists an ordinal $\xi'$ such that $\lambda > \xi' \geq \xi$ and $F(\xi') = \alpha'$, and thus we can assume that $F(\xi) = \alpha'$. The commutative square
\[\begin{tikzcd}
	{\mathrm{dom}(u_{\alpha'} \mathrel{\Box} g)} & {D\xi} \\
	{V_{\alpha'} \odot Y} & L
	\arrow["{w'}", from=1-1, to=1-2]
	\arrow["{u_{\alpha'} \mathrel{\Box} g}"', from=1-1, to=2-1]
	\arrow["{\varphi_\xi}", from=1-2, to=2-2]
	\arrow["{(p_2 \cdot v)^*}"', from=2-1, to=2-2]
\end{tikzcd}\]
is one of the squares from the construction, and thus there exists $\beta < \gamma_\xi$ such that $(g,(p_2 \cdot v)^*,w') = (g_\beta, v_\beta, w_\beta)$. Define the diagonal that we are looking for to be the following composite.
\[\begin{tikzcd}
	{V_{\alpha'} \odot Y} & {D_{\beta + 1}(\xi + 1)} & {D(\xi + 1)} & {\operatorname{colim}_{\alpha < \lambda} D\alpha}
	\arrow["{h'_\beta}", from=1-1, to=1-2]
	\arrow["{\iota'_{\beta + 1}}", from=1-2, to=1-3]
	\arrow["{\iota_{\xi + 1}}", from=1-3, to=1-4]
\end{tikzcd}\]
This diagonal indeed makes the two triangles commute, since
$$m \cdot \iota_{\xi + 1} \cdot \iota'_{\beta + 1} \cdot h'_{\beta} = \varphi_{\xi + 1} \cdot \iota'_{\beta + 1} \cdot h'_\beta = \varphi_{\xi + 1, \beta + 1} \cdot h'_{\beta} = v_\beta = (p_2 \cdot v)^*,$$
and
\begin{align*}
\iota_{\xi + 1} \cdot \iota'_{\beta + 1} \cdot h'_{\beta} \cdot (u_{\alpha'} \mathrel{\Box} g) &= \iota_{\xi + 1} \cdot \iota'_{\beta + 1} \cdot f'_{\beta, \beta + 1} \cdot f'_{0,\beta} \cdot w_\beta \\
&= \iota_{\xi + 1} \cdot \iota'_{\beta} \cdot f'_{0,\beta} \cdot w_\beta \\
&= \iota_{\xi + 1} \cdot \iota'_0 \cdot w_\beta \\
&= \iota_{\xi} \cdot w_\beta \\
&= (w^*,(p_1 \cdot v)^*).
\end{align*}
Thus, the proof of $m \in \mathcal{I}^{\overset{\mathscr{F}}{\pitchfork}}$ is finished.
\epf
\rem
Note that $\mathcal{I} \subseteq {}^{\overset{\mathscr{F}}{\pitchfork}}(\mathcal{I}^{\overset{\mathscr{F}}{\pitchfork}})$ and $\mathcal{I} \subseteq (^{\overset{\mathscr{F}}{\pitchfork}}\mathcal{I})^{\overset{\mathscr{F}}{\pitchfork}}$. In fact, this holds for any binary relation and $\overset{\mathscr{F}}{\pitchfork}$ is a binary relation. Therefore, from the stability properties in Section \ref{sec:stabilityproperties} and the fact that $\mathcal{I} \subseteq {}^{\overset{\mathscr{F}}{\pitchfork}}(\mathcal{I}^{\overset{\mathscr{F}}{\pitchfork}})$, we obtain an enriched weak $\mathscr{F}$-factorization system $(^{\overset{\mathscr{F}}{\pitchfork}}(\mathcal{I}^{\overset{\mathscr{F}}{\pitchfork}}), \mathcal{I}^{\overset{\mathscr{F}}{\pitchfork}})$ by Theorem \ref{thm:soa}.
\endrem
\rem\label{rem:replacesmallnessbypres}
The smallness assumption in Theorem \ref{thm:soa} holds in particular if
\begin{enumerate}[(i)]
\item the domains and codomains of morphisms in $\mathcal{I}$ are presentable in the enriched sense, and
\item the domains and codomains of morphisms in $\mathcal{J}$ are presentable in the unenriched sense.
\end{enumerate}
\endrem
\rem\label{rem:restrictcolim}
By inspecting the construction in the proof of Theorem \ref{thm:soa} we see that instead of assuming the existence of all pushouts and transfinite composites in $\mathcal{K}$ it suffices to assume that $\mathcal{K}$ has the following special types of pushouts and transfinite composites in order to conclude that Theorem \ref{thm:soa} holds:
\begin{enumerate}[(i)]
\item pushouts of the form $\mathrm{dom}(u \mathrel{\Box} g)$ from Definition \ref{def:uodotf}, where $u \in \mathcal{J}$, $g \in \mathcal{I}$,
\item pushouts of morphisms of the form $u \mathrel{\Box} g$, where $u \in \mathcal{J}$, $g \in \mathcal{I}$, and
\item transfinite composites of morphisms from $(\mathcal{I} \mathrel{\mathrm{rel}} \mathcal{J})\text{-}\mathrm{cell}$.
\end{enumerate}
\endrem
\examples\label{ex:enrichedsoafinal} In all of the following examples (see Examples \ref{ex:liftings}.(\ref{ex:liftingsSet})-(\ref{ex:liftingsVIso}) for the corresponding $\mathscr{F}$-lifting properties) Theorem \ref{thm:soa} holds. In Examples (\ref{ex:insSet})-(\ref{ex:insSSet}), $\mathcal{V}$ is indeed a cosmos, and in Example (\ref{ex:insPres}) it is assumed to be. In all examples we will assume that $\mathscr{F} = ({}^\pitchfork(\mathcal{J}^\pitchfork), \mathcal{J}^\pitchfork)$.
\begin{enumerate}[(1)]
\item\label{ex:insSet} $\mathcal{V} = \mathbf{Set}$, $\mathcal{J} = \{u \colon \emptyset \to 1\}$. The stability of $\mathscr{L}$ under corners follows from the fact that $u \mathrel{\Box} u = u$ and from Lemma \ref{lem:strnabla}. Thus, from Theorem \ref{thm:soa} we obtain the classical small object argument \cite[Theorem 2.1.14]{Hov99} for weak factorization systems, since the domains of morphisms in $\mathcal{I}$ being small relative to $\mathcal{I}\text{-}\mathrm{cell}$ \cite[Definition 2.1.9]{Hov99} corresponds precisely to the smallness assumption in Theorem \ref{thm:soa}. Note that instead of assuming the existence of all colimits in $\mathcal{K}$ it suffices to assume that $\mathcal{K}$ has pushouts of morphisms in $\mathcal{I}$ and transfinite composites of morphisms in $\mathcal{I}\text{-}\mathrm{cell}$, see Remark \ref{rem:restrictcolim}.\par
If we instead choose ${\mathcal{J} = \{u \colon \emptyset \to 1, v \colon 2 \to 1\}}$, then the stability of $\mathscr{L}$ under corners follows from the fact that $\mathscr{L}$ is the class of all functions, since $\mathscr{R}$ is the class of all bijections. If $\mathcal{K}$ is locally presentable, then by Theorem \ref{thm:soa} we obtain a variant of the small object argument that gives us an orthogonal factorization system $(^{\bot}(\mathcal{I}^{\bot}), \mathcal{I}^{\bot})$ (where $\bot$ denotes the corresponding $\mathscr{F}$\nobreakdash-lifting property), which forms part of \mbox{\cite[Theorem 2.2]{FR08}}. Note that in order to apply Theorem \ref{thm:soa} to obtain such an orthogonal factorization system it suffices to assume, instead of local presentability of $\mathcal{K}$, that $\mathcal{K}$ has certain special types of pushouts, coequalizers, and transfinite composites (Remark \ref{rem:restrictcolim}), and that the domains and codomains of morphisms in $\mathcal{I}$ are small relative to $(\mathcal{I} \mathrel{\mathrm{rel}} \mathcal{J})\text{-}\mathrm{cell}$.
\item\label{ex:insCat} $\mathcal{V} = \mathbf{Cat}$, $\mathcal{J} = \{u \colon \emptyset \to 1, v \colon 2 \to \mathbf{2}, w \colon \mathbf{2}' \to \mathbf{2}\}$. Note that $\emptyset$, $1$, $2$, $\mathbf{2}$, $\mathbf{2'}$ are finitely presentable as objects of $\mathbf{Cat}$, $\mathcal{J}$ has cardinality less than $\aleph_0$, and $\mathscr{L}$ is stable under corners. The stability of $\mathscr{L}$ under corners follows from \cite[Theorem 5.1]{Rez96}.\par
\item\label{ex:insGrpd} $\mathcal{V} = \mathbf{Grpd}$, $\mathcal{J} = \{u \colon \emptyset \to 1, v \colon 2 \to \mathbf{2}_g, w \colon \mathbf{2}'_g \to \mathbf{2}_g\}$. Note that $\emptyset$, $1$, $2$, $\mathbf{2}_g$, $\mathbf{2}'_g$ are finitely presentable as objects of $\mathbf{Grpd}$, $\mathcal{J}$ has cardinality less than $\aleph_0$, and $\mathscr{L}$ is stable under corners. The stability of $\mathscr{L}$ under corners can be shown completely analogously as in Example (\ref{ex:insCat}).\par
The small object argument for (2,1)-categories that we obtain in this way differs from the the small object argument for (2,1)-categories in \cite[Section 5]{Kan21}, since the $\mathscr{F}$-lifting property that we are considering has a 2-dimensional aspect that's not present in the aforementioned paper. Even if we remove the previously mentioned 2\nobreakdash-dimensional aspect by omitting $v$ and $w$ from $\mathcal{J}$, then the arguments are still different, since the one in the aforementioned paper has homotopical aspects that are not present in our argument.
\item\label{ex:insCh} $\mathcal{V} = \mathbf{Ch}$, $\mathcal{J} = \{S^{n - 1} \hookrightarrow D^n \mid n \in \mathbb{Z}\}$. Note that for each $n \in \mathbb{Z}$, $S^{n -1}$ and $D^n$ are finitely presentable (and hence $\aleph_1$-presentable) as objects of $\mathbf{Ch}$, and $\mathcal{J}$ has cardinality less than $\aleph_1$. The stability of the class $\mathscr{L}$ under corners follows from \cite[Proposition 4.2.13]{Hov99}.
\item\label{ex:insSSet} $\mathcal{V} = \mathbf{SSet}$, $\mathcal{J} = \{\partial \Delta^n \hookrightarrow \Delta^n \mid n \geq 0\}$. Note that for each $n \geq 0$, $\partial \Delta^n$ and $\Delta^n$ are finitely presentable (and hence $\aleph_1$-presentable) as objects of $\mathbf{SSet}$, and $\mathcal{J}$ has cardinality less than $\aleph_1$. The stability of the class $\mathscr{L}$ under corners follows from \cite[Proposition 4.2.8]{Hov99}.
\item\label{ex:insPres} $\mathcal{V}_0$ is a locally $\lambda$-presentable category with a set $\mathcal{V}_\lambda$ of $\lambda$-presentable objects that form a strong generator, where $\lambda$ is a regular cardinal. Consider $$\mathcal{J} = \{u_V \colon \emptyset \to V \mid V \in \mathcal{V}_\lambda\} \cup \{\nabla_V \colon V + V \to V \mid V \in \mathcal{V}_\lambda\}.$$ Let $\mu$ be a regular cardinal greater than $\max\{\lambda, |\mathcal{J}|\}$. Then all the domains and codomains of morphisms in $\mathcal{J}$ are $\mu$-presentable and $\mathcal{J}$ has cardinality less than $\mu$. The stability under corners follows from the fact that $\mathscr{L}$ is the class of all morphisms in $\mathcal{V}$.
\end{enumerate}
\endexamples
\example
In every example in Examples \ref{ex:enrichedsoafinal} the morphism $!_I \colon \emptyset \to I$ belongs to $\mathscr{L}$, which is the condition that appears in Proposition \ref{prop:selfenrichment}. However, Theorem \ref{thm:soa} doesn't require $!_I \in \mathscr{L}$ to hold. We mention a simple example that demonstrates this: Let $\mathcal{V} = \mathbf{Set}$ and ${\mathcal{J} = \{v \colon 2 \to 1\}}$, i.e.\ $\mathscr{F} = (\text{surjective}, \text{injective})$. Note that $2$ and $1$ are finitely presentable in $\mathbf{Set}$ and $\mathcal{J}$ has cardinality less than $\aleph_0$. The stability of $\mathscr{L}$ under corners follows from Lemma \ref{lem:strnabla} and the fact that $\mathscr{L}$ is the class of all surjections. Indeed, the codomain of $v \mathrel{\Box} v$ is $1$ and the domain of $v \mathrel{\Box} v$ is non-empty, since it is a pushout of non-empty sets. We also remark that $f \overset{\mathscr{F}}{\pitchfork} k$ if and only if for each commutative $(f,k)$\nobreakdash-square (\ref{cmd:fksquare}) it is the case that if $d,d' \colon B \to C$ are two diagonals making the triangles inside the square commute, then $d = d'$. Note that existence of the diagonal is not required.
\endexample
\example
Every example in Examples \ref{ex:enrichedsoafinal} involves a locally presentable cosmos $\mathcal{V}$. In some cases it is also possible to apply the Theorem \ref{thm:soa} when $\mathcal{V}$ is not locally presentable. We illustrate this by the following example: Let $\mathcal{V} = \mathbf{Top}$ be the cosmos of compactly generated topological spaces with the standard cartesian closed monoidal structure and let $\mathcal{K} = [\mathbf{2},\mathbf{Top}]$ be the $\mathbf{Top}$-category of arrows in $\mathbf{Top}$. Furthermore, let $\mathcal{J} = \{i_n \colon S^{n - 1} \hookrightarrow D^n \mid n \geq 0\}$ be the standard generating cofibrations on $\mathbf{Top}$ and let $\mathcal{I} = \{\mathrm{id}_{i_0} \colon i_0 \to i_0, !_{i_0} \colon \mathrm{id}_\emptyset \to i_0\}$. We now explain what $\mathcal{I}$ consists of in terms of commutative squares. We denote the empty topological space by $\emptyset$ and the one-point topological space by $1$. Recall that $S^{-1} = \emptyset$ and $D^0 = 1$, i.e.\ $i_0$ is the unique morphism $\emptyset \to 1$. Thus $\mathcal{I}$ consists of the following two commutative squares.
\[\begin{tikzcd}
	\emptyset & \emptyset && \emptyset & \emptyset \\
	1 & 1 && \emptyset & 1
	\arrow["{\mathrm{id}_{\emptyset}}", from=1-1, to=1-2]
	\arrow["{i_0}"', from=1-1, to=2-1]
	\arrow["{i_0}", from=1-2, to=2-2]
	\arrow["{\mathrm{id}_\emptyset}", from=1-4, to=1-5]
	\arrow["{\mathrm{id}_\emptyset}"', from=1-4, to=2-4]
	\arrow["{i_0}", from=1-5, to=2-5]
	\arrow["{\mathrm{id}_1}"', from=2-1, to=2-2]
	\arrow["{i_0}"', from=2-4, to=2-5]
\end{tikzcd}\]
Then $\mathscr{L}$ is stable under corners \cite[Proposition 4.2.11]{Hov99}. Moreover, every topological space is small relative to closed $T_1$-inclusions \cite[Lemma 2.4.1]{Hov99}, and thus every object of $[\mathbf{2},\mathbf{Top}]_0$ is small relative to componentwise closed $T_1$-inclusions. Note that ${i_n \mathrel{\Box} \mathrm{id}_{i_0} = \mathrm{id}_{i_n}}$ and $i_n \mathrel{\Box} \operatorname{!}_{i_0} = \operatorname{!}_{i_n}$. Therefore, each morphism in $(\mathcal{I} \mathrel{\mathrm{rel}} \mathcal{J})\text{-}\mathrm{cell}$ is a componentwise closed $T_1$-inclusion according to \cite[Lemma 2.4.5]{Hov99}. Thus, we can apply Theorem \ref{thm:soa} to obtain the small object argument for the projective model structure on $[\mathbf{2},\mathbf{Top}]$ from \cite[Theorem 3.1]{Hov14}.
\endexample
\begin{references*}

\bibitem[AR94]{AR94}
Ji\v{r}\'{\i} Ad\'{a}mek and Ji\v{r}\'{\i} Rosick\'{y}.
\newblock {\em Locally presentable and accessible categories}, volume 189 of
  {\em London Mathematical Society Lecture Note Series}.
\newblock Cambridge University Press, Cambridge, 1994.
\newblock \href {https://doi.org/10.1017/CBO9780511600579}
  {\path{doi:10.1017/CBO9780511600579}}.

\bibitem[Day70]{Day70}
Brian Day.
\newblock On closed categories of functors.
\newblock In {\em Reports of the {M}idwest {C}ategory {S}eminar, {IV}}, Lecture
  Notes in Math., Vol. 137, pages 1--38. Springer, Berlin-New York, 1970.
\newblock \href {https://doi.org/10.1007/BFb0060438}
  {\path{doi:10.1007/BFb0060438}}.

\bibitem[Day74]{Day74}
Brian Day.
\newblock On adjoint-functor factorisation.
\newblock In {\em Category {S}eminar ({P}roc. {S}em., {S}ydney, 1972/1973)},
  Lecture Notes in Math., Vol. 420, pages 1--19. Springer, Berlin-New York,
  1974.
\newblock \href {https://doi.org/10.1007/BFb0063097}
  {\path{doi:10.1007/BFb0063097}}.
  
\bibitem[FR08]{FR08}
L.~Fajstrup and J.~Rosick\'{y}.
\newblock A convenient category for directed homotopy.
\newblock {\em Theory Appl. Categ.}, 21:No. 1, 7--20, 2008.

\bibitem[Gar09]{Gar09}
Richard Garner.
\newblock Understanding the small object argument.
\newblock {\em Appl. Categ. Structures}, 17(3):247--285, 2009.
\newblock \href {https://doi.org/10.1007/s10485-008-9137-4}
  {\path{doi:10.1007/s10485-008-9137-4}}.

\bibitem[Hov99]{Hov99}
Mark Hovey.
\newblock {\em Model categories}, volume~63 of {\em Mathematical Surveys and
  Monographs}.
\newblock American Mathematical Society, Providence, RI, 1999.
\newblock \href {https://doi.org/10.1090/surv/063}
  {\path{doi:10.1090/surv/063}}.
  
\bibitem[Hov14]{Hov14}
Mark Hovey.
\newblock Smith ideals of structured ring spectra, 2014.
\newblock \href {https://arxiv.org/abs/1401.2850} {\path{arXiv:1401.2850}}.

\bibitem[JK01]{JK01}
G.~Janelidze and G.~M. Kelly.
\newblock A note on actions of a monoidal category.
\newblock {\em Theory Appl. Categ.}, 9:61--91, 2001.

\bibitem[Kan21]{Kan21}
Kristóf Kanalas.
\newblock The (2,1)-category of small coherent categories, 2021.
\newblock \href {https://arxiv.org/abs/2104.13239} {\path{arXiv:2104.13239}}.

\bibitem[Kel82]{Kel82}
G.~M. Kelly.
\newblock Structures defined by finite limits in the enriched context. {I}.
\newblock {\em Cahiers Topologie G\'{e}om. Diff\'{e}rentielle}, 23(1):3--42,
  1982.
\newblock Third Colloquium on Categories, Part VI (Amiens, 1980).

\bibitem[Kel05]{Kel05}
G.~M. Kelly.
\newblock Basic concepts of enriched category theory.
\newblock {\em Repr. Theory Appl. Categ.}, (10):vi+137, 2005.
\newblock Reprint of the 1982 original [Cambridge Univ. Press, Cambridge;
  MR0651714].

\bibitem[LW14]{LW14}
Rory B.~B. Lucyshyn-Wright.
\newblock Enriched factorization systems.
\newblock {\em Theory Appl. Categ.}, 29:No. 18, 475--495, 2014.

\bibitem[MMSS01]{MMSS01}
M.~A. Mandell, J.~P. May, S.~Schwede, and B.~Shipley.
\newblock Model categories of diagram spectra.
\newblock {\em Proc. London Math. Soc. (3)}, 82(2):441--512, 2001.
\newblock \href {https://doi.org/10.1112/S0024611501012692}
  {\path{doi:10.1112/S0024611501012692}}.

\bibitem[MU22]{MU22}
Dylan McDermott and Tarmo Uustalu.
\newblock Flexibly graded monads and graded algebras.
\newblock In {\em Mathematics of program construction}, volume 13544 of {\em
  Lecture Notes in Comput. Sci.}, pages 102--128. Springer, Cham, 2022.
\newblock \href {https://doi.org/10.1007/978-3-031-16912-0_4}
  {\path{doi:10.1007/978-3-031-16912-0_4}}.

\bibitem[Qui67]{Qui67}
Daniel~G. Quillen.
\newblock {\em Homotopical algebra}.
\newblock Lecture Notes in Mathematics, No. 43. Springer-Verlag, Berlin-New
  York, 1967.
\newblock \href {https://doi.org/10.1007/BFb0097438}
  {\path{doi:10.1007/BFb0097438}}.

\bibitem[Rez96]{Rez96}
Charles Rezk.
\newblock A model category for categories, 1996.
\newblock
  \url{https://ncatlab.org/nlab/files/Rezk_ModelCategoryForCategories.pdf}.
  Accessed: 11-16-2023.

\bibitem[Rie14]{Rie14}
Emily Riehl.
\newblock {\em Categorical homotopy theory}, volume~24 of {\em New Mathematical
  Monographs}.
\newblock Cambridge University Press, Cambridge, 2014.
\newblock \href {https://doi.org/10.1017/CBO9781107261457}
  {\path{doi:10.1017/CBO9781107261457}}.

\bibitem[RT07]{RT07}
Ji\v{r}\'{\i} Rosick\'{y} and Walter Tholen.
\newblock Factorization, fibration and torsion.
\newblock {\em J. Homotopy Relat. Struct.}, 2(2):295--314, 2007.

\bibitem[Sub21]{Sub21}
Chaitanya~Leena Subramaniam.
\newblock From dependent type theory to higher algebraic structures.
\newblock Ph.D. thesis. Université de Paris, 2021.
\newblock \href {https://arxiv.org/abs/2110.02804} {\path{arXiv:2110.02804}}.

\bibitem[Woo76]{Woo76}
Richard~James Wood.
\newblock Indicial methods for relative categories.
\newblock Ph.D. thesis. Dalhousie University, 1976.
\newblock \url{https://dalspace.library.dal.ca//handle/10222/55465}.

\end{references*}

\end{document}